\magnification 1150 \baselineskip 12 pt

\def \b {\beta}

\def \la {\lambda}
\def \a {\alpha}

\def \g {\gamma}

\input psfig.sty

\def \sect#1{\bigskip  \noindent{\bf #1} \medskip }
\def \subsect#1{\bigskip \noindent{\it #1} \medskip}
\def \th#1#2{\medskip \noindent {\bf Theorem #1.}   \it #2 \rm \medskip}
\def \prop#1#2{\medskip \noindent {\bf Proposition #1.}   \it #2 \rm \medskip}
\def \cor#1#2{\medskip \noindent {\bf Corollary #1.}   \it #2 \rm \medskip}
\def \pf {\noindent  {\it Proof}.\quad }
\def \lem#1#2{\medskip \noindent {\bf Lemma #1.}   \it #2 \rm \medskip}

\def\sqr#1#2{{\vcenter{\vbox{\hrule height.#2pt\hbox{\vrule width.#2pt height#1pt \kern#1pt\vrule width.#2pt}\hrule height.#2pt}}}}

\def \square{\hfill\mathchoice\sqr56\sqr56\sqr{4.1}5\sqr{3.5}5}

\headline={\ifnum\pageno=1 \hfill \else \hfill {\rm \folio} \fi}

\centerline{\bf Pricing Options in Incomplete Equity Markets}
\medskip \centerline{\bf via the Instantaneous Sharpe Ratio}
\bigskip
\centerline{Erhan Bayraktar}

\centerline{Virginia R. Young} \bigskip

\centerline{Department of Mathematics, University of Michigan}
\centerline{Ann Arbor, Michigan, 48109} \bigskip

\centerline{Version: 10 June 2007} \bigskip

\noindent{\bf Abstract:}  We use a continuous version of the standard deviation premium principle for pricing in incomplete equity markets by assuming that the investor issuing an unhedgeable derivative security requires compensation for this risk in the form of a pre-specified instantaneous Sharpe ratio.  First, we apply our method to price options on non-traded assets for which there is a traded asset that is correlated to the non-traded asset. Our main contribution to this particular problem is to show that our seller/buyer prices are the upper/lower good deal bounds of Cochrane and Sa\'{a}-Requejo (2000) and of Bj\"{o}rk and Slinko (2006) and to determine the analytical properties of these prices. Second, we apply our method to price options in the presence of stochastic volatility.  Our main contribution to this problem is to show that the instantaneous Sharpe ratio, an integral ingredient in our methodology, is the negative of the market price of volatility risk, as defined in Fouque, Papanicolaou, and Sircar (2000).

\bigskip

\noindent{\it Keywords:} Pricing derivative securities, incomplete markets, Sharpe ratio, correlated assets, stochastic volatility, non-linear partial differential equations, good deal bounds.

\sect{1. Introduction}

In this paper, we provide new insights for pricing and hedging derivative securities in incomplete equity markets. We consider two sources of incompleteness: (1) the asset underlying the derivative security is non-traded, as in Davis (2000) and Musiela and Zariphopoulou (2004), or (2) the underlying asset exhibits stochastic volatility.  In the first of these two problems, an investor partially hedges her position on a derivative written on the non-traded asset by trading a correlated asset. Executive options and weather derivatives are examples of options on non-traded assets.  An investor holding an executive option is usually not allowed to hedge her position using the underlying asset (see e.g.\ Leung and Sircar (2006)).  In the case of a temperature derivative, there is no underlying asset available, and an agent trading in such a derivative might want to hedge her position by trading in electricity futures, which are correlated to the temperature.  Also, in many situations, the underlying asset can be traded, but an investor might choose to hedge the option on the underlying by using an index because of high transaction costs or illiquidity, for example. A real option is also a derivative on a non-traded asset.

In the second of these problems, an investor attempts to hedge against the stochastic volatility in the underlying asset. Stochastic volatility models are proposed in the literature to account for the implied volatility smile; see Fouque, Papanicolaou, and Sircar (2000) for a review.

For each of these two problems, we determine the price/value $P=\{P_t\}_{t \geq 0}$ of the option by assuming that the investor minimizes the local variance of a suitably defined portfolio, $\Pi=\{\Pi_t\}_{t \geq 0}$, and is compensated for the residual risk by specifying a so-called {\it instantaneous Sharpe ratio} of this portfolio. When the investor is a market maker and can set prices, $P$ will indeed be the price of the derivative. In this paper, we assume that the seller determines the price process $P$ so that the (instantaneous) Sharpe ratio of a portfolio that is composed of the option liability and a self-financing strategy that minimizes the local variance of the portfolio is equal to a pre-specified constant $\alpha$, which can be thought of as the risk loading the seller charges for taking on extra risk that she cannot hedge.  Let us explain this further: To each self-financing strategy $V$ and price process $P$, we associate a portfolio process $\Pi$ (of the option liability and the self-financing strategy). We only consider pairs $(V,P)$ such that the corresponding portfolio has instantaneous Sharpe ratio $\alpha$. Among all these pairs $(V,P)$ we choose the one whose corresponding portfolio has the smallest instantaneous variance; that is, we choose a pair $(V,P)$ so that the portfolio $\Pi$ is on the efficient frontier. Therefore, we not only find the price (process) of the option, but also its hedge (see Section 2.2).

Our pricing mechanism is similar to the way insurance companies set prices for their products, namely by charging the expected value of the payoff plus an additional amount for the risk.  In our setting, some of the risk the option writer undertakes can be hedged by the market, and for the unhedgeable risk the writer of the option takes the attitude of an insurance company.  The parameter $\alpha$ can be chosen to be a function of the underlying equities, but this would not change the analysis or the qualitative conclusions significantly.  If the investor cannot set the prices, then we view $P$ as a value that reflects her risk preferences because $P$ is what she would charge/pay if she could sell/buy the derivative given her risk preferences, which is reflected in $\alpha$.

$P$ can be thought of as the price the investor thinks is fair.  Once $P$ is determined, she can use it to make buying or selling decisions given the market price of the derivative. In this sense, the value process $P$ is similar to the utility indifference price (see e.g.\ Musiela and Zariphopoulou (2004); Sircar and Zariphopoulou (2005)): both of these prices reflect the attitude of the investor towards risk, and investors choose these prices as benchmarks when making investment decisions.  Pricing/valuation using the instantaneous Sharpe ratio might prove to be more user-friendly (that is, understandable and implementable by traders) than the utility indifference price, since the risk preferences in the former are specified in terms of a widely known Sharpe ratio, whereas in the latter the risk preferences are specified in terms of a utility function (which is usually chosen to be exponential so that the price is time consistent, see e.g.\ Cheridito and Kupper (2006)).  An investor may not be able to determine her risk aversion precisely, but anybody who is familiar with trading knows her risk-return tradeoff.

The price/value satisfies a non-linear partial differential equation (PDE). In Section 2, we derive the conditions under which this equation has a unique solution. In our framework, for any European option, the seller's price is strictly greater than the buyer's price (see Corollary 2.8 and Section 2.5); hence, there is a positive bid-ask spread, as in utility indifference pricing (see e.g.\ Musiela and Zariphopoulou (2004); Sircar and Zariphopoulou (2005)). However, unlike in indifference pricing (under exponential utility, the buyer/seller price is strictly concave/convex in the number of options bought/sold, see e.g.\ Ilhan et al.\ (2004)), our pricing PDE satisfies the scaling property (Proposition 2.1), so that the price is linear in terms of number of units traded, which is consistent with market prices.  The two properties of our prices, the linearity of the prices in the number of units bought or sold and the bid-ask spread, are shared by the prices obtained using the (static) standard deviation premium principle (in fact our pricing mechanism is a dynamic version of this principle), in which  the seller/buyer price $P(X)$ of a random variable $X$ is given by
$
P(X)={\bf E}(X) \, \pm \, \alpha \sqrt{{\bf Var}(X)}
$ (see e.g.\ Gerber (1979)), and its variant, the financial standard deviation principle (see Schweizer (2001), pages 44-45).

Our approach to pricing in incomplete markets is much like pricing under the minimal martingale measure, which also minimizes the local variance of a hedged portfolio (see e.g.\ F\"ollmer and Schweizer (1991)). The new twist in our work is that the writer or buyer of the option receives an additional compensation through a pre-specified Sharpe ratio $\alpha$. Indeed, when $\alpha=0$ and the investor is risk neutral, $P^{\alpha0}$ is equal to the expected value of payoff of the option under the minimal martingale measure; see (2.23). When $\alpha>0$, then the buyer's price and the seller's price satisfy $P^{b}<P^{\alpha0}<P^{s}$; see Corollary 2.8 and (2.31). The parameter $\alpha$ can be thought of as the parameter of risk aversion; see Theorem 2.7.

By using comparison principles for PDEs, as found in Barles et al.\ (2003) and Walter (1970), in Section 2, we show that the price we derive for an option on a non-traded underlying satisfies a number of appealing properties. In the special case for which the prices of both the traded and the non-traded assets follow geometric Brownian motion, we give a closed form solution for the price of a European put option written on the non-traded asset in terms of the Black-Scholes formula; see Corollary 2.6. This clearly distinguishes our pricing mechanism from the utility indifference pricing, since in the latter framework explicit expressions for the prices of such options are not available.

Windcliff et al.\ (2007) and Forsyth and Labahn (2006) also analyze the special case for which the prices of both the traded and the non-traded assets follow geometric Brownian motion. They independently derive the same pricing PDE, which has only one state variable besides the time variable in this special case, and present a numerical algorithm to solve it. They provide numerical algorithms to solve both European and American options in this context.  While most of the focus of Windcliff et al.\ (2007) is on the numerical algorithm to solve the non-linear pricing PDE, our paper explores analytical properties of the option prices: (1) We derive the pricing PDE under very general modeling assumptions and derive conditions under which the non-linear pricing PDE has a unique solution (Theorem A.1).  (2) Under the general modeling assumptions, we show that there is positive bid-ask spread (see Corollary 2.8 and Section 2.5) although the buyer/seller price scale linearly in the number of the options sold (bought) (Proposition 2.1).  (3) We show that discounted expected value of the option pay-off under the minimal martingale measure is always between the buyer's and the seller's prices for any value of $\a$.  (4) By using a comparison principle for non-linear PDEs (we provide the conditions under which the comparison principle holds), we show how the option contract behaves with respect to the model parameters mostly when the Delta of the contract never changes sign.  (5) Lastly, we provide an explicit expression for the price of put and call contracts when both assets are correlated geometric Brownian motions (Corollaries 2.8 and 2.16).  An interesting observation is that there is a parity between the buyer's call price and the seller's put price, and the buyer's put price and the seller's call price.

Although our pricing mechanism is very different--we use a continuous version of the standard deviation premium principle--it turns out that our buyer/seller price corresponds to the lower/upper good deal bound of Cochrane and Sa\'{a}-Requejo (2000) and of Bj\"{o}rk and Slinko (2006).  The no-arbitrage price interval of a future pay-off $X$ is $(\inf_{{\bf Q} \in {\cal M}} E^{Q}(e^{-rT} X)$, $\sup_{{\bf Q} \in \cal M} {\bf E}^{Q}(e^{-rT}X))$, in which ${\cal M}$ is the set of martingale measures, is a wide interval, and Cochrane and Sa\'{a}-Requejo (2000) find a reasonably small interval for the prices by ruling out the extremely good deals by putting a bound on the absolute value of the volatility/market price of risk (with respect to the independent Brownian motion driving the non-traded asset) of the Radon-Nikodym derivative of the measures ${\bf Q} \in {\cal M}$ with respect to ${\bf P}$. The lower and upper prices for the resulting pricing interval are obtained through  solving stochastic control problems, although this is not precisely stated in Cochrane and Sa\'{a}-Requejo (2000). Our pricing mechanism provides a new interpretation of the lower and upper prices of Cochrane and Sa\'{a}-Requejo (2000) as the buyer (bid) and seller (ask) prices under the risk aversion parameter $\a$. Our main contribution is determining the analytical properties of the prices as we mentioned in the previous paragraph.  Also, see the closing remark of Section 2.5.

Our other contribution is exhibited in Section 3, where we consider the case for which the option liability is written on a tradable asset that exhibits stochastic volatility.  We show that the price of the option for an investor with a given risk tolerance (that is, the instantaneous Sharpe ratio in our context) coincides with the price under a risk-neutral measure with a market price of volatility risk that is equal to the investor's risk tolerance (see the remark on pages 20-21 following (3.16)).  Our framework provides an explanation for the market price of volatility risk and, therefore, the implied volatility smile in terms of traders' risk preferences. This means that the lower and upper good deal bounds are attained by the martingale measures whose market price of risk is equal to the instantaneous Sharpe ratio of the buyer and the seller, respectively.

The rest of the paper is organized as follows: In Section 2, we apply our method to price a European option on an asset that is not traded but is correlated to an asset that can be traded.  In Section 2.1, we introduce the market model, and in Section 2.2, we derive the PDE for the seller's price.  In the remainder of Section 2, we give conditions under which the pricing PDE has a unique solution, calculate the prices of European options explicitly in some specific cases, and look at the comparative statistics of the price. We show that there is a positive bid-ask spread, and the pricing mechanism has several other appealing properties. In Section 3, we apply our method to price a European option on a traded asset that exhibits stochastic volatility and show that specifying a Sharpe ratio in our pricing mechanism is the same as specifying the market price of risk. Section 4 concludes
the paper.

\sect{2. Pricing Derivative Securities on a Non-Traded Asset}

In this section, we present the financial market for the writer of an option on a non-traded asset.  We obtain the hedging strategy for the writer of the option when the writer can trade on a correlated asset.  We describe how to use the instantaneous Sharpe ratio to price the option and derive the resulting partial differential equation that the price solves.  By using comparison arguments, we determine many properties of the price.  We also obtain the corresponding equation that the buyer's price solves and examine the bid-ask spread.

\subsect{2.1. Financial Market}

Suppose a writer issues an option on an non-traded asset whose price process $S$ follows
$$
dS_t = \mu(S_t, t) \, S_t \, dt + \sigma(S_t, t) \, S_t \, dZ_t, \eqno(2.1)$$

\noindent in which $Z$ is a standard Brownian motion on a probability space $(\Omega, {\cal F}, {\bf P})$.   The European option pays $g(S_T)$ at time $T$.  Because the writer cannot trade on the asset underlying the derivative security, the writer trades on a correlated asset whose price process $H$ follows
$$
dH_t = a(H_t, t) \, H_t \, dt + b(H_t, t) \, H_t \, dW_t, \eqno(2.2)$$

\noindent in which $W$ is a standard Brownian motion on $(\Omega, {\cal F}, {\bf P})$ with constant coefficient of correlation $\rho$ with respect to the Brownian motion $Z$. 
We assume that the coefficients $\mu$, $\sigma$, $a$, and $b$ satisfy growth conditions and are locally Lipschitz in $S$ or $H$, as appropriate.  These conditions ensure the existence and uniqueness of the solutions of the two SDEs above.

\subsect{2.2. Writer's Price $($Ask Price$)$}

The writer faces the unhedgeable risk of not being able to trade on the security underlying the payout $g(S_T)$; therefore, the writer demands a return greater than the return $r$ on the riskless asset.  One measure of the risk that the writer assumes is the local standard deviation of a suitably-defined portfolio. A natural tie between the excess return and the standard deviation is the ratio of the former to the latter, the so-called {\it instantaneous Sharpe ratio}.  In what follows, we find the hedging strategy to minimize the local variance of the  portfolio, then we set the price of the option so that the resulting instantaneous Sharpe ratio is equal to a given value.

Denote the value (price) of the option by $P = P(S, H, t)$, in which we explicitly recognize that the price of the option depends on the non-traded asset's price $S$ and the traded asset's price $H$ at time $t$.  Suppose the writer creates a portfolio $\Pi$ with value $\Pi_t$ at time $t$.  The portfolio contains the obligation to pay $g(S_T)$ at time $T$.  Additionally, the writer holds a self-financing portfolio $V=\{V_t\}_{t \geq 0}$ composed of shares in the traded asset and a money market account.  Let $r \geq 0$ be the interest rate earned by the money market
account, and let $\pi_t$ denote the number of shares of the traded asset that the investor holds at time $t \geq 0$. The dynamics of the self-financing portfolio $V$ are given by
$$
dV_t=\left( r V_t+ (a(H_t,t)-r) \, \pi_t \, H_t \right)dt+ b(H_t,t) \, \pi_t \, H_t \, dW_t, \quad t \geq 0. \eqno(2.3)$$

By applying It\^o's Lemma (see e.g.\ Karatzas and Shreve (1991)) to the portfolio process  $\Pi_t = -P(S_t, H_t, t) + V_t$, $t \geq 0$, we obtain
$$
\eqalign{\Pi_{t+h} &= \Pi_t - \int_t^{t+h} {\cal D}^{\mu, a} P(S_s, H_s, s) \, ds + \int_t^{t+h} b(H_s, s) H_s (\pi_s - P_H(S_s, H_s, s)) \, dW_s \cr
& \quad - \int_t^{t+h} \sigma(S_s, s) S_s P_S(S_s, H_s, s) \, dZ_s + \int_t^{t+h} \pi_s (a(H_s, s)-r) H_s \, ds, + \int_t^{t+h}r \Pi_s ds,}  \eqno(2.4)$$

\noindent in which ${\cal D}^{f, k}$, with $f = f(S, t)$ and $k = k(H, t)$ deterministic functions, is an operator defined on the set of appropriately differentiable functions on $G = {\bf R}^+ \times {\bf R}^+ \times [0, T]$ by
$$
{\cal D}^{f, k} v = v_t + f S v_S + k H v_H + {1 \over 2} \sigma^2 S^2 v_{SS} + \rho \sigma b SH v_{SH} + {1 \over 2} b^2 H^2 v_{HH}-r v.  \eqno(2.5)$$

We next calculate the expectation and variance of $\Pi_{t+h}$ conditional on the information available at time $t$, namely ${\cal F}_t$.  First, given $\Pi_t = \Pi$, define the stochastic process $Y_h$ for $h \ge 0$ by
$$
Y_h = \Pi - \int_t^{t+h} {\cal D}^{\mu, a} P(S_s, H_s, s) ds + \int_t^{t+h} [\pi_s (a(H_s, s)-r) H_s+r \Pi_s] \, ds. \eqno(2.6)$$

\noindent  Thus, ${\bf E}(\Pi_{t+h} | {\cal F}_t) = {\bf E}^{S, H, t} (Y_h)$, in which ${\bf E}^{S, H, t}$ denotes the conditional expectation given $S_t = S$ and $H_t = H$.  From (2.4), we have
$$
\eqalign{\Pi_{t + h} &= Y_h + \int_t^{t+h} b(H_s, s) H_s (\pi_s - P_H(S_s, H_s, s)) \, dW_s \cr & \quad - \int_t^{t+h} \sigma(S_s, s) S_s P_S(S_s, H_s, s) \, dZ_s.} \eqno(2.7)$$

\noindent It follows that
$$
\eqalign{& {\bf Var}(\Pi_{t+h} | {\cal F}_t) = {\bf E}((\Pi_{t+h} - {\bf E}Y_h)^2 | {\cal F}_t) \cr
& \quad = {\bf E}^{S, H, t} (Y_h - {\bf E}Y_h)^2 + {\bf E}^{S, H, t} \int_t^{t+h} b^2(H_s, s) H_s^2 (\pi_s - P_H(S_s, H_s, s))^2 ds \cr &
\qquad - 2 \rho {\bf E}^{S, H, t} \int_t^{t+h} \sigma(S_s, s) b(H_s, s) S_s H_s P_S(S_s, H_s, s) (\pi_s - P_H(S_s, H_s, s)) ds \cr
& \qquad + {\bf E}^{S, H, t}  \int_t^{t+h} \sigma^2(S_s, s) S_s^2 P_S^2(S_s, H_s, s) ds +o(h).} \eqno(2.8)$$

\noindent We choose $\pi_t$ in order to minimize the local variance $\lim_{h \rightarrow 0} {1 \over h} {\bf Var}(\Pi_{t+h} | {\cal F}_t)$, a measure of risk of the portfolio; therefore, the optimal investment in the traded asset at time $t$ is
$$
\pi^*_t = P_H(S_t, H_t, t) + \rho \, {\sigma(S_t, t) \over b(H_t, t)} \, {S_t \over H_t} \, P_S(S_t, H_t, t). \eqno(2.9)$$

\noindent Under this assignment, the drift and local variance become, respectively,
$$
\lim_{h \rightarrow 0} {1 \over h} ({\bf E}(\Pi_{t+h} | {\cal F}_t) - \Pi) = - {\cal D}^{\mu - \rho \sigma (a-r)/b, r} P(S, H,
t)+r \Pi , \eqno(2.10)$$

\noindent and
$$
\lim_{h \rightarrow 0} {1 \over h} {\bf Var}(\Pi_{t+h} | {\cal F}_t) =  (1 - \rho^2) \sigma^2(S, t) S^2 P_S^2(S, H, t).
\eqno(2.11)$$

Now, we come to pricing via the instantaneous Sharpe ratio. Because the minimum variance is positive, the writer is unable to completely hedge the risk of the option written on the non-traded security.  Therefore, the price should reimburse the writer for this risk, say, by a constant multiple $\a$ of the local standard deviation of the portfolio.  It is this $\a$ that is the instantaneous Sharpe ratio.

From (2.11), we learn that the local standard deviation of the portfolio equals
$$
\lim_{h \rightarrow 0} \sqrt{{1 \over h}  {\bf Var}(\Pi_{t+h} | {\cal F}_t)} = \sqrt{1 - \rho^2} \, \sigma(S, t) \, S \, \big|
P_S(S, H, t) \big|.  \eqno(2.12)$$

\noindent To determine the value (price) $P$, we set the drift of the portfolio $\Pi$ equal to the short rate times the portfolio {\it plus} $\a$ times the local standard deviation.  Thus, from (2.10) and (2.12), we have that $P$ solves the equation
$$
- {\cal D}^{\mu - \rho \sigma (a-r)/b, r} P(S, H, t)+r \Pi  = r \Pi + \a \sqrt{1 - \rho^2} \, \sigma(S, t) \, S \, \big| P_S(S, H,
t) \big|, \eqno(2.13)$$

\noindent for a given $\a \ge 0$.  It follows that the writer's price $P = P(S, H, t)$ solves the non-linear PDE given by
$$
\left\{ \eqalign{&P_t + \left( \mu(S, t) - (a(H, t) - r) \rho {\sigma(S, t) \over b(H, t)} \right) S P_S + r H P_H \cr &
\quad + {1 \over 2} \sigma^2(S, t) S^2 P_{SS} + \rho \sigma(S, t) b(H, t) S H P_{SH} + {1 \over 2} b^2(H, t) H^2 P_{HH} - rP \cr &
\quad = -\a \sqrt{1 - \rho^2} \, \sigma(S, t) \, S \, \big| P_S \big| , \cr
& P(S, H, T) = g(S).} \right. \eqno(2.14)$$

One can think of the right-hand side of the PDE in (2.14) as adding a margin to the return of the portfolio because the risk arising from the non-traded asset is not completely hedgeable.  If there were no risk loading, that is, if $\a = 0$, then the price is such that the expected return on the price is $r$.  If $\a > 0$, then the expected return on the price is greater than $r$.  Therefore, $\a$, the Sharpe ratio, measures the degree to which the writer's total expected return is in excess of $r$, as a proportion of the standard deviation of the return. On the other hand, when the assets are perfectly correlated, that is, $|\rho|=1$, then the price does not carry any risk loading, and as a result of Feynman-Kac Theorem, the price is an expectation under the risk-neutral measure. Therefore, our framework produces the correct price when the option contract is perfectly hedgeable.

\subsect{2.3. Qualitative Properties of the Risk-Adjusted Price}

In this section, we discuss qualitative properties of the risk-adjusted price $P$ in (2.14).  To begin, we have the
following proposition whose proof is clear, so we omit it.

\prop{2.1}{Suppose $P^c$ is the price, as determined by the method in Section 2.2, for an option with payment $c \, g(S_T)$ at time $T$, with $c \ge 0$.  Then, $P^c = c \, P$, in which $P$ solves $(2.14)$.}

\medskip

From Proposition 2.1, we learn that the price for an option scales by the number of units of the option.  Therefore, if one were to sell two options each paying $g(S_T)$ at time $T$, then the price would be twice that of the single option.  This scaling property does not hold for utility indifference prices (Musiela and Zariphopoulou, 2004).

In what follows, we show that we can determine the sign of $P_S$ in some cases so that equation (2.14) becomes a linear equation, and we examine how the price $P$ responds to changes in the model parameters.  To this end, we need a comparison principle, and in this section, we rely on a comparison principle from Barles et al.\ (2003); see Appendix A.

\medskip

\noindent{\bf Assumption 2.2.} Henceforth, in Section 2, we assume that:
\item{(1)} $\sigma$ and $b$ satisfy [i] and [ii] in Theorem A.1 and are differentiable with respect to their first variable.
\item{(2)} $\tilde{\mu}$ satisfies the  conditions in the hypothesis of Lemma A.2 and is differentiable with respect to its first and second variable.
\item{(3)} A generic payoff $g$ satisfies the following growth condition:  There exist constants $L > 0$ and $\la \ge 1$ such that
$$
g(S) \le L \left(1+( \ln S)^{2 \lambda}\right) $$ for all $S > 0$.
\item{(4)}  The functions $\mu$, $\sigma$, $a$, and $b$ are bounded.
\item{(5)} The PDE in (2.14) satisfies a uniform ellipticity condition, that is, there is a constant $\delta > 0$ such that for any $(\xi_{1},\xi_2) \in {\bf R}^2$
$$
\sigma^2(e^{x},t) \xi^2_1+ 2 \rho \sigma(e^x,t) b(e^y,t) \xi_1 \xi_2+ b^{2}(e^y,t) \xi_2^2 \ge \delta \sqrt{\xi_1^2+\xi_2^2}.$$

\medskip

We now apply Theorem A.1 and Lemma A.2 repeatedly to determine qualitative properties of the price $P$.  In our first application, we show that if $g$ is monotone, then $P$ is monotone in $S$. In this case, (2.14) reduces to a linear PDE.

\th{2.3}  {Assume that $\mu$ and $\sigma$ do not depend on $S$ and that the payoff $g$ is continuously differentiable. If $g'(S) \ge (\le) \; 0$ on ${\bf R}^+,$ then $P_S \ge (\le) \; 0$ on $G$.}

\pf We outline the proof of this assertion. First note that under Assumption 2.2,  the PDE in (2.14) has a unique viscosity solution in $\cal{C}$.  Suppose that the payoff function $g$ is increasing with respect to $S$, and consider the solution $f = f(S, H, t)$ of
$$
\left\{ \eqalign{&f_t + \left( \tilde{\mu}(H, t) + \a \sqrt{1 - \rho^2} \, \sigma(t) \right) S f_S + r H f_H \cr
& \quad + {1 \over 2} \sigma^2(t) S^2 f_{SS} + \rho \sigma(t) b(H, t) S H f_{SH} + {1 \over 2} b^2(H, t) H^2 f_{HH} - rf = 0,\cr & f(S, H, T) = g(S),} \right. \eqno(2.15)$$

\noindent in which $\tilde \mu(H, t) = \mu(t) - (a(H, t) - r) \rho \, \sigma(t)/ b(H, t)$.

Under Assumption 2.2, it follows from Friedman (1975, Theorem 5.3) that (2.15) has a unique solution among functions that satisfy a certain growth condition. Moreover, this solution is $C^{2, 2,1}({\bf R}^+ \times {\bf R}^+ \times [0,T])$ and its derivative satisfies the same growth condition.

By differentiating the PDE in (2.15) with respect to $S$, we get a linear homogeneous PDE for the derivative of $f$, which we denote by $F = f_S$. The terminal condition for the PDE that $F$ satisfies is $F(S, H, T) = g'(S)$. Under Assumption 2.2, it follows from Theorem A.1 that $F$'s PDE has a unique continuous viscosity solution that satisfies the growth condition in Theorem A.1 [v].  (Because $f_S$ satisfies the growth condition, it is the viscosity solution of the PDE that we derived by differentiation.)

By defining a differential operator $\cal L$ via $F$'s PDE, we get that ${\cal L} F = {\cal L} {\bf 0} = 0$ on $G$, in which $\bf 0$ denotes the function that is identically 0.  Thus, by applying the comparison result in Theorem A.1, we conclude that $F = f_S \ge 0$ on $G$.

Recall that Theorem A.1 and Lemma A.2 imply that the solution of (2.14) is unique.  It follows that $f = P$ and $P_S \ge 0$ on $G$.  We obtain a parallel result when $g' \le 0$, by replacing $\a$ with $-\a$ in (2.15).   $\square$

\medskip

\noindent {\bf Remark:} (1)  The assumptions that $g$ is continuously differentiable and that $g$ satisfies Assumption 2.2(3) are sufficient but not necessary conditions for Theorem 2.3 to hold.  The reason we assume that $g$ is continuously differentiable is so that we can apply the comparison principle (or a maximum principle) for parabolic differential equations to the PDE that the derivative of the solution of (2.15), namely $f_S$, satisfies, which requires that the terminal condition $g'(S)$ be continuous. However, we can derive conclusions beyond what the maximum principle tells us:

\item{(i)} In the case of a put option, $g$ is not differentiable. However, in Corollary 2.5, by using the fact that put payouts can be uniformly approximated by smooth functions (from above), we show that the statement of Theorem 2.3 is still valid.

\item{(ii)} In the case of a call option, again $g$ is not differentiable; moreover, it does not satisfy Assumption 2.2(3). However, In Corollary 2.14 below, by using a parity relationship between the {\it buyer's} put option price and the seller's call option price we show that the statement of Theorem 2.3 still holds.

(2) Also, the assumption that $\mu$ and $\sigma$ do not depend on $S$ is a sufficient but not a necessary condition for the results in Theorem 2.3 to hold. This assumption is a sufficient condition to ensure that the derivative of the solution of (2.15) solves the PDE obtained by differentiating (2.15).  By other means one can prove that this is true without assuming $\mu$ and $\sigma$ do not depend on $S$. Indeed, sufficient conditions for the results of Theorem 2.3 to hold are the following two growth conditions: $|\sigma_S| S \le K(1 + \min[\ln S, (\ln S)^2])$ and $|\tilde{\mu}_S| S \le K(1 + (\ln  S)^2)$. This observation follows from Walter (1970, Section 28, pages 213-215). In Corollary 2.4, we assume that the hypotheses of Theorem 2.3 hold, but one could instead assume that Walter's conditions hold and obtain the same result.

\medskip

From Theorem 2.3, we have the intuitive result that if the payoff $g$ is increasing (decreasing) with the price $S$ of the underlying security, then the price $P$ of the derivative security is also increasing (decreasing) with $S$.  One can also show that if $g$ is monotone and if $g'' \ge (\le) \; 0$, then under appropriate growth conditions on the coefficients of the PDE that $P_{SS}$ satisfies (so that the comparison principle can be applied), $P_{SS} \ge
(\le) \; 0$.

In the following corollary, we apply the Feynman-Kac Theorem to represent the price $P$ when $g' \ge 0$.

\cor{2.4} {Suppose that the hypotheses of Theorem $2.3$ hold, and suppose that $g' \ge 0$ on ${\bf R}^+$.  Then, we have
$$
P(S, H, t) = {\bf \hat E}^{S, H, t} [e^{-r(T - t)} g(S_T)], \eqno(2.16)$$
in which $S$ and $H$ follow the processes
$$
dS_s = \left(\tilde{\mu}(H_s, s) + \a \sqrt{1 - \rho^2} \, \sigma(s) \right) S_s \, ds + \sigma(s) \, S_s \, d \hat Z_s, \eqno(2.17)$$
with $\hat Z_s = Z_s + \int_0^s  {\mu(u) - \tilde{\mu}(H_u, u) \over \sigma(u)} \, du - \a \sqrt{1 - \rho^2} \; s$, and
$$
dH_s = r \, H_s \, ds + b(H_s, s) \, H_s \, d \hat W_s, \eqno(2.18)$$
with $\hat W_s = W_s + \int_0^s {a(H_u, u) - r \over b(H_u, u)} \, du$.  The processes $\hat Z$ and $\hat W$ are standard Brownian motions on the probability space $(\Omega, {\cal F}, {\bf \hat P})$, in which the Radon-Nikodym derivative of $\bf \hat P$ with respect to $\bf P$ is
$$
\eqalign{{d {\bf \hat P} \over d {\bf P}} \bigg|_{{\cal F}_t} &= \exp \left( - \int_0^t \left( {\mu(s) - \tilde{\mu}(H_s, s) \over \sigma(s)} - \a \sqrt{1 - \rho^2} \right) \, dZ_s  \right) \cr
& \quad \times \exp \left( - {1 \over 2} \int_0^t \left( {\mu(s) - \tilde{\mu}(H_s, s) \over \sigma(s)} - \a \sqrt{1 - \rho^2} \right)^2 \, ds \right) \cr
& \quad \times \exp \left( - \int_0^t {a(H_s, s) - r \over b(H_s, s)} \, dW_s - {1 \over 2} \int_0^t \left( {a(H_s, s) - r \over b(H_s, s)} \right)^2 \, ds  \right),} \eqno(2.19)$$
and ${\bf \hat E}$ denotes expectation with respect to $\bf \hat P$.}

\pf Because $P_S \ge 0$ on $G$, (2.14) becomes a linear PDE.  The result, then, follows by applying the Feynman-Kac Theorem to this linear PDE; see e.g.\ Karatzas and Shreve (1991, pages 366-368). $\square$

\medskip

We have a representation parallel to (2.16) when $g' \le 0$ by replacing $\a$ with $-\a$ in (2.17), which we do not state for the sake of brevity.

\cor{2.5} {Suppose that $g$ is the payoff for a European put option with strike price $K$, that is, $g(S) = (K - S)_+$.  Then, we can represent the solution of $(2.14)$ as in $(2.16)$ with $\a$ replaced by $-\a$ in $(2.17)$. Moreover, the put option price satisfies $P_{S} \leq 0$.}

\pf We cannot directly apply Corollary 2.4 to represent $P$ because $g(S) = (K - S)_+$ is not differentiable.  We approximate $g$ with smooth functions and apply a limit theorem to prove our assertion.  Indeed, define $g_n$, for $n > 1/K$, by
$$
g_n(S) = \cases{K - S, &if  $S < K - {1 \over n}$, \cr {n \over 4} \left( K - {1 \over n} - S \right)^2, &if $K - {1 \over n} \le S < K + {1 \over n}$, \cr 0, &if $S \ge K + {1 \over n}$.} $$

\noindent Then, $g_n$ is smooth, and $\lim_{n \rightarrow \infty} g_n(S) = (K - S)_+$.  It follows from the Lebesgue Dominated Convergence Theorem (Royden, 1968) that
$$
\lim_{n \rightarrow \infty} {\bf \hat E}^{S, H, t} [e^{-r(T - t)} g_n(S_T)] = {\bf \hat E}^{S, H, t} [e^{-r(T - t)} (K - S_T)_+], \eqno(2.20)$$

\noindent in which $S$ follows the process in (2.17) with $\a$ replaced by $-\a$, and $H$ follows the process in (2.18).

Because ${\partial \over \partial S} {\bf \hat E}^{S, H, t} \left[e^{-r(T - t)} g_n(S_T)\right] \le 0$, it follows that the limit $P_S \le 0$.  Thus, $P$ solves (2.15) with $\a$ replaced by $-\a$ and with $P_S \le 0$, from which it follows that $P$ solves (2.14) by uniqueness of the solution of (2.14).  In other words, the price for a European put option can be represented as in (2.16) with $\a$ replaced by $-\a$ in (2.17).  $\square$

\medskip

One can also show that for $g(S) = (K - S)_+$, then $P_{SS} \ge 0$, as is true for the ordinary Black-Scholes price of a  put option.

When $\mu$, $\sigma$, $a$, and $b$ are constant, then from Corollaries 2.4 and 2.5, one can derive an explicit expression for the price of a European put option by using the Black-Scholes pricing formula for a put option written on a security that pays dividends.  Note that the price is independent of $H$ if $a$ and $b$ are
independent of $H$; we revisit this fact in Section 2.4.  Also, note that when all the coefficients are constants, Assumption 2.2 is trivially satisfied.

\cor{2.6} {When $\mu,$ $\sigma,$ $a,$ and $b$ are constant $($that is, when the both the underlying asset and the traded asset are geometric Brownian motions$),$ then the price $P$ of a European put option is given by
$$
P(S, t) = K  e^{-r \tau} \Phi \left({\ln \left({K \over S} \right) - (r - \delta - {1 \over 2} \sigma^2) \tau \over \sigma \sqrt{\tau}}\right) - S e^{-\delta \tau} \Phi \left({\ln \left({K \over S }\right) - (r-\delta + {1 \over 2} \sigma^2) \tau \over \sigma \sqrt{\tau}}\right), \eqno(2.21)$$
in which $\tau = T - t$, $\delta = r - \tilde{\mu} + \alpha \sqrt{1-\rho^2} \, \sigma,$ and $\Phi$ denotes the cumulative distribution function of the standard normal random variable. }
\medskip

\noindent {\bf Remark:} Note that when $\rho = 1$, the market becomes complete and in order to avoid arbitrage one imposes that ${a - r \over b} = {\mu - r \over \sigma}$; see equation (2.3) in Davis (2000).  Thus, when $\rho = 1$ and ${a - r \over b} = {\mu - r \over \sigma}$, then $\delta = 0$ and (2.21) reduces to the ordinary Black-Scholes price for a European put option.  Similarly, if $\rho = -1$, then the necessary condition to ensure that the market is free of arbitrage is ${a - r \over b} = - {\mu - r \over \sigma}$; again, $\delta = 0$ and (2.21) reduces to the Black-Scholes price for a put option.
\medskip

\subsect{2.4. Variations with respect to the Model Parameters}

Our next results show that as we vary the model parameters, the price $P$ responds consistently with what we expect.  Unless stated otherwise, $\mu$ and $\sigma$ are functions of $S$ and $t$, and $a$ and $b$ are functions of $H$ and $t$.

\th{2.7} {Suppose $0 \le \a_1 <  \a_2 $, and let $P^{\a_i}$ be the solution to $(2.14)$ with $\a = \a_i,$ for $i = 1, 2$. Then, $P^{\a_1} \le P^{\a_2}$ on $G$.}

\pf  Define a differential operator $\cal L$ on $\cal G$ by (A.1) and (A.8) with $\a = \a_1$.  Because $P^{\a_1}$ solves (2.14) with $\a = \a_1$, we have ${\cal L} P^{\a_1} = 0$.  Also,
$$
{\cal L} P^{\a_2} = -(\a_2 - \a_1) \sqrt{1 - \rho^2} \, \sigma S \big| P_S \big| \le 0 = {\cal L} P^{\a_1}. \eqno(2.22)$$

\noindent In addition, both $P^{\a_1}$ and $P^{\a_2}$ satisfy the terminal condition $P^{\a_i}(S, H, T) = g(S)$.  Thus, Theorem A.1 and Lemma A.2 imply that $P^{\a_1} \le P^{\a_2}$ on $G$. $\square$

\medskip

Theorem 2.7 states that as the parameter $\a$ increases, the price $P$ increases.  For this reason, we refer to $P$ as the {\it risk-adjusted price}.  We have the following corollary to Theorem 2.7.

\cor{2.8} {Let $P^{\a0}$ be the solution to $(2.14)$ with $\a = 0;$ then, $P^{\a0} \le P^{\a}$ for all $\a \ge 0$, and we can express the lower bound $P^{\a0}$ as follows:
$$
P^{\a0}(S, H, t) = {\bf \hat E}^{S, H, t} [e^{-r(T - t)} g(S_T)], \eqno(2.23)$$
\noindent in which $S$ and $H$ follow the processes given in $(2.16)$ and $(2.17)$, respectively, with $\a = 0$ in $(2.16)$.}

\medskip

\pf Theorem 2.7 implies that $P^{\a0} \le P^{\a}$ for all $\a \ge 0$, and by substituting $\a = 0$ in (2.14), the Feynman-Kac Theorem (Karatzas and Shreve, 1991) implies the representation of $P^{\a0}$ in (2.23).  $\square$

\medskip

Corollary 2.8 justifies calling $P$ a risk-adjusted price because the lower bound is attained when $\a = 0$. One might call this lower bound, the expected value of the pay-off under the minimal martingale measure (see F\"{o}llmer and Schweizer (1991)), a {\it risk-neutral} price, since the option writer is not charging anything for the unhedged risk exposure. Corollary 2.8 also justifies the use of the phrase {\it risk parameter} when referring to $\a$.  We can think of $P^\a - P^{\a0}$ as the {\it risk charge} to compensate the writer of this option for the unhedgeable risk of writing an option on a non-traded security.

Next, we examine how the risk-adjusted price $P$ varies with the drifts of $S$ and $H$.

\th{2.9} {Suppose $r \le \mu_1(S, t) \le \mu_2(S, t)$ on ${\bf R}^+ \times [0, T]$, and let $P^{\mu_i}$ denote the solution to $(2.14)$ with $\mu = \mu_i,$  for $i = 1, 2$.  If $P^{\mu_i}_S \ge (\le) \hbox{ } 0$ on $G$ for $i = 1$ or $2,$ then $P^{\mu_1} \le (\ge) \; P^{\mu_2}$ on $G$.}

\pf  Suppose $P_S^{\mu_2} \ge 0$ on $G$.  Define a differential operator $\cal L$ on $\cal G$ by (A.1) and (A.8) with $\mu = \mu_1$.  Because $P^{\mu_1}$ solves (2.14) with $\mu = \mu_1$, we have ${\cal L} P^{\mu_1} = 0$.  Also, because $P_S^{\mu_2} \ge 0$, we have
$$
{\cal L} P^{\mu_2} = -(\mu_2 - \mu_1) S P^{\mu_2}_S \le 0 = {\cal L} P^{\mu_1}. \eqno(2.24)$$

\noindent In addition, both $P^{\mu_1}$ and $P^{\mu_2}$ satisfy the terminal condition $P^{\mu_i}(S, H, T) = g(S)$. Theorem A.1 and Lemma A.2 imply that $P^{\mu_1} \le P^{\mu_2}$ on $G$.  The other cases follow similarly.  $\square$

\medskip

Theorem 2.9 is an intuitively pleasing result:  If $\mu_1 \le \mu_2$, then $\mu_2$ will make $S$ increase more than $\mu_1$ will.  Moreover, if $P$ responds positively to changes in $S$, then $P^{\mu_1} \le P^{\mu_2}$ is a natural conclusion.

\th{2.10} {Suppose $r \le a_1(H, t) \le a_2(H, t)$ on ${\bf R}^+ \times [0, T]$, and let $P^{a_i}$ denote the solution to $(2.14)$ with $a = a_i,$  for $i = 1, 2$.  If $\rho P_S^{a_i} \ge (\le) \hbox{ } 0$ on $G$ for $i = 1$ or $2,$ then $P^{a_1} \ge (\le) \; P^{a_2}$ on $G$.}

\pf Suppose $\rho P_S^{a_2} \ge 0$ on $G$.  Define a differential operator $\cal L$ on $\cal G$ by (A.1) and (A.8) with $a = a_1$.  Because $P^{a_1}$ solves (2.14) with $a = a_1$, we have ${\cal L} P^{a_1} = 0$.  Also, because $\rho P_S^{a_2} \ge 0$, we have
$$
{\cal L} P^{a_2} = (a_2 - a_1) \rho {\sigma \over b} S P_S^{a_2} \ge 0 = {\cal L} P^{a_1}. \eqno(2.25)$$

\noindent In addition, both $P^{a_1}$ and $P^{a_2}$ satisfy the terminal condition $P^{a_i}(S, H, T) = g(S)$. Theorem A.1 and Lemma A.2 imply that $P^{a_1} \ge P^{a_2}$ on $G$.  The other cases follow similarly.  $\square$

\medskip

Theorem 2.10 tells us that if $S$ and $H$ are positively correlated, and if $P_S \ge 0$, then the price decreases with respect to the drift on the hedging asset.

\medskip

\noindent 2.4.1.  Further Properties of $P$ when $a$ and $b$ are Independent of $H$

\medskip

If the drift and volatility of $H$ are independent of $H$, that is, if $a(H, t) = a(t)$ and $b(H, t) = b(t)$, then the price $P$ is independent of $H$.  Indeed, suppose $P(S, H, t) = P(S, t)$, and substitute that expression into (2.14) to obtain
$$
\left\{ \eqalign{&P_t + \left( \mu(S, t) - (a(t) - r) \rho {\sigma(S, t) \over b(t)} \right) S P_S + {1 \over 2} \sigma^2(S,
t) S^2 P_{SS} - rP  \cr
& \quad = - \a \sqrt{1 - \rho^2} \, \sigma(S, t) \, S \, \big| P_S \big| , \cr
& P(S, T) = g(S).} \right. \eqno(2.26)$$

\noindent Note that (2.26) is independent of $H$.  Thus, by the uniqueness of the solutions of (2.14) and (2.26),  the solutions of those two PDEs are equal.  We now examine how $P$ varies with respect to changes in the volatility of $S$ and $H$.

\th{2.11} {Suppose $0 \le \sigma_1(S, t) \le \sigma_2(S, t)$ on ${\bf R}^+ \times [0, T]$, and let $P^{\sigma_i}$ denote the solution to $(2.26)$ with $\sigma = \sigma_i,$ for $i = 1, 2$.  If $P^{\sigma_i}_S \ge 0$ and $P^{\sigma_i}_{SS} \ge 0$ for $i = 1$ or $2,$ and if $\a \sqrt{1 - \rho^2} \ge  \rho (a - r)/b,$ then $P^{\sigma_1} \le P^{\sigma_2}$ on $G$.}

\pf Suppose that $P^{\sigma_2}_S \ge 0$ and $P^{\sigma_2}_{SS} \ge 0$.  Define a differential operator $\cal L$ on $\cal G$ by (A.1) and (A.8) with $\sigma = \sigma_1$.  Because $P^{\sigma_1}$ solves (2.26) with $\sigma = \sigma_1$, we have ${\cal L} P^{\sigma_1} = 0$.  Also, because $\a \sqrt{1 - \rho^2} \ge  \rho (a - r)/b$, we have
$$
\eqalign{{\cal L} P^{\sigma_2} &= - {1 \over 2} (\sigma_2^2 - \sigma_1^2) S^2 P_{SS}^{\sigma_2} - (\sigma_2 - \sigma_1) (\a \sqrt{1 - \rho^2} - \rho (a - r)/b) S P^{\sigma_2}_S \cr
& \le 0 = {\cal L} P^{\sigma_1}.} \eqno(2.27)$$

\noindent In addition, both $P^{\sigma_1}$ and $P^{\sigma_2}$ satisfy the terminal condition $P^{\sigma_i}(S, T) = g(S)$.  Theorem A.1 and Lemma A.2 imply that $P^{\sigma_1} \le P^{\sigma_2}$ on $G$.  The case for which $P^{\sigma_1}_S \ge 0$ and $P^{\sigma_1}_{SS} \ge 0$ follows similarly.  $\square$

\medskip

Note that the condition $\a \sqrt{1 - \rho^2} \ge  \rho (a - r)/b$ in Theorem 2.11 is automatic if $\rho \le 0$ because we assume that $a \ge r$.  Otherwise, the condition holds if the instantaneous Sharpe ratio $\a$ (or level of risk aversion) is large enough relative to the Sharpe ratio corresponding to $H$'s price process.

\th{2.12} {Suppose $r \le b_1(H, t) \le b_2(H, t)$ on ${\bf R}^+ \times [0, T]$, and let $P^{b_i}$ denote the solution to $(2.26)$ with $b = b_i,$  for $i = 1, 2$.  If $\rho P_S^{b_i} \ge (\le) \hbox{ } 0$ on $G$ for $i = 1$ or $2,$ then $P^{b_1} \le (\ge) \; P^{b_2}$ on $G$.}

\pf  Suppose $\rho P_S^{b_2} \ge 0$ on $G$.  Define a differential operator $\cal L$ on $\cal G$ by (A.1) and (A.8) with $b = b_1$.  Because $P^{b_1}$ solves (2.26) with $b = b_1$, we have ${\cal L} P^{b_1} = 0$.  Also, because $\rho P_S^{b_2} \ge 0$, we have
$$
{\cal L} P^{b_2} = - (a - r) \rho \left( {\sigma \over b_1} - {\sigma \over b_2} \right) S P_S^{b_2} \le 0 = {\cal L} P^{b_1}. \eqno(2.28)$$

\noindent In addition, both $P^{b_1}$ and $P^{b_2}$ satisfy the terminal condition $P^{b_i}(S, T) = g(S)$. Theorem A.1 and Lemma A.2 imply that $P^{b_1} \le P^{b_2}$ on $G$.  The other cases follow similarly.  $\square$

\medskip

Theorem 2.12 tells us that if $S$ and $H$ are positively correlated, and if $P_S \ge 0$, then the price increases with the volatility $b$ of the hedging asset.  To end this section, we determine how $P$ varies with the correlation coefficient $\rho$.

\th{2.13} {Suppose $0 \le \rho_1 \le \rho_2$, and let $P^{\rho_i}$ denote the solution to $(2.26)$ with $\rho = \rho_i,$  for $i = 1, 2$. If $P^{\rho_i}_S \ge 0$ for $i = 1$ or $2,$ then $P^{\rho_1} \ge P^{\rho_2}$ on $G$.}

\pf  Suppose $P_S^{\rho_2} \ge 0$ on $G$.  Define a differential operator $\cal L$ on $\cal G$ by (A.1) and (A.8) with $\rho = \rho_1$.  Because $P^{\rho_1}$ solves (2.26) with $\rho = \rho_1$, we have ${\cal L} P^{\rho_1} = 0$.  Also, we have
$$
{\cal L} P^{\rho_2} = \left( (\rho_2 - \rho_1) (a - r) {\sigma \over b} + \a \left( \sqrt{1 - \rho_1^2} - \sqrt{1 - \rho_2^2}
\right) \right) S P_S^{\rho_2} \ge 0 = {\cal L} P^{\rho_1}. \eqno(2.29)$$

\noindent In addition, both $P^{\rho_1}$ and $P^{\rho_2}$ satisfy the terminal condition $P^{\rho_i}(S, T) = g(S)$. Theorem A.1 and Lemma A.2 imply that $P^{\rho_1} \ge P^{\rho_2}$ on $G$.  The case for which $P_S^{\rho_1} \ge 0$ follows similarly.  $\square$

\medskip

Theorem 2.13 is intuitively pleasing.  If $P_S \ge 0$, then as the correlation coefficient increases, the traded asset becomes a better hedge for the derivative security written on the non-traded asset; therefore, the price decreases.  We have a similar result for $\rho_2 \le \rho_1 \le 0$, which we omit for the sake of brevity.

\subsect{2.5.  Buyer's Price $($Bid Price$)$ $P^b$}

In this section, we consider the problem from a buyer's point of view and determine a price, $P^b$, which we call the buyer's price. In this problem, we consider an investor who will receive a payment of $g(S_T)$ at time $T$. Similar to the set up of the seller's problem, the buyer holds a self-financing portfolio that is composed of shares in the traded asset $H$ and a money market account to minimize the local variance of her entire portfolio. As the seller does, the buyer also specifies her risk aversion through an instantaneous Sharpe ratio, $\beta$, and considers a price fair for the option, if her portfolio's instantaneous Sharpe ratio is equal to $\beta$.  Again, among the portfolios with instantaneous Sharpe ratio equal to $\beta$, the investor chooses  the portfolio with the minimum local variance.  To determine the buyer's price $P^b$, one can follow the argument in Section 2.2 to obtain that $P^b$ solves
$$
\left\{ \eqalign{&P^b_t + \left( \mu(S, t) - (a(H, t) - r) \rho {\sigma(S, t) \over b(H, t)} \right) S P^b_S + r H P^b_H \cr
& \quad + {1 \over 2} \sigma^2(S, t) S^2 P^b_{SS} + \rho \sigma(S, t) b(H, t) S H P^b_{SH} + {1 \over 2} b^2(H, t) H^2 P^b_{HH} - rP^b \cr
& \quad = \beta \sqrt{1 - \rho^2} \, \sigma(S, t) \, S \, \big| P^b_S \big| \cr & P^b(S, H, T) = g(S).} \right. \eqno(2.30)$$

\noindent In other words, $P^b$ solves (2.14) with $\a \ge 0$ replaced by $-\beta$.  Alternatively, one can observe that $-P^b$ solves (2.14) with $g$ replaced by $-g$ and again reach the conclusion that $P^b$ solves (2.30).  The results of the previous subsections modify easily to apply to the buyer's price, so we do not repeat them here.

The most interesting difference between the writer's price $P$ and the buyer's price $P^b$ is that $P^b$ {\it decreases} with respect to $\beta$, while $P$ {\it increases} with respect to $\a$.  Thus, for any $\beta \ge 0$ and $\a \ge 0$, we have that

$$P^{b,\beta} \le P^{b,\a0} = P^{\a0} \le P^\a, \eqno(2.31)$$

\noindent in which $P^{\a0}$ is the price when $\a = 0 = \beta$.  Thus, the bid-ask spread $P^\a - P^{b,\beta}$ is always non-negative, as one expects of a reasonable pricing model. Here we denoted the buyer's price $P^b$ as $P^{b,\beta}$ to emphasize its dependence on the instantaneous Sharpe ratio $\beta$.

\medskip

\noindent{\bf Remark:} Note that the buyer's price for a European put option, when $\mu$, $\sigma$, $a$, and $b$ are constant, equals the expression in (2.21) with $\alpha$ replaced by $-\beta$.

\cor{2.14} {When $\mu,$ $\sigma,$ $a,$ and $b$ are constant, the seller's price for a European call option is given by
$$
C^s(S, t) = S e^{-\delta \tau} \Phi \left({\ln \left({ S \over K}\right) +(r-\delta + {1 \over 2} \sigma^2) \tau \over \sigma \sqrt{\tau}}\right)-K  e^{-r \tau} \Phi \left({\ln \left({ S \over K} \right) + (r -\delta - {1 \over 2} \sigma^2) \tau \over \sigma \sqrt{\tau}}\right),  \eqno(2.32)$$
\noindent in which $\tau = T - t$, $\delta = r -\tilde{\mu} - \alpha \sqrt{1-\rho^2} \, \sigma.$}

\pf  Denote the solution of (2.15) by ${\tilde C}^s$ when $g(S) = (S - K)_+$. Note that the terms in (2.15)  that have differential with respect to $H$ disappear.  The function ${\tilde C}^s$ satisfies a parity relationship with the {\it buyer's} put price, $P^b$, when $\beta=\alpha$, namely
$$
{\tilde C}^s + Ke^{-r \tau} = P^b + S e^{-\delta \tau},
$$
as can be see from their respective PDEs.  From this parity relationship, we see that $\tilde{C}^b$ is equal to the right-hand side of (2.32) and that ${\tilde C}^s_S=P^b_S+e^{- \delta \tau}$.  By using the remark right before the statement of this corollary, we can directly differentiate $P^b$ and obtain that $P^b_S \geq -e^{-\delta \tau}$ from which it follows that ${\tilde C}^{s}_S \geq 0$. Therefore ${\tilde C}^{s}$ is a solution of (2.14) with $g(S)=(S-K)_+$ and the comparison result leads to the uniqueness, from which it follows that ${\tilde C}^s = C^s$, since the latter is the unique solution of (2.14) with $g(S)=(S-K)_+$. $\square$

\medskip

\noindent {\bf Remark:}  To find the buyer's call option price, when her risk aversion is measured by $\beta$, simply replace $\a$ by $-\beta$ in the expression in (2.32).

\medskip

\noindent{\bf Numerical Example:} Figure 1 shows how the price of a European call option changes with respect to some of the parameter values.  The purpose of this figure is to illustrate the magnitude of the effects of the parameters in addition to confirming the results of Section 2.4.  We use the same parameters as in Table 1 of Windcliff et al.\ (2007), namely $S=100, K=100, T=1, r=0.05, \sigma=0.2, a=0.077$, and $b=0.3$. For the surface plot of European call option values, we set $mu=0.07$; for the second plot, we set $\rho=0.9$ and $\alpha=0.2$.

\centerline{\vbox{\hbox{\psfig{figure=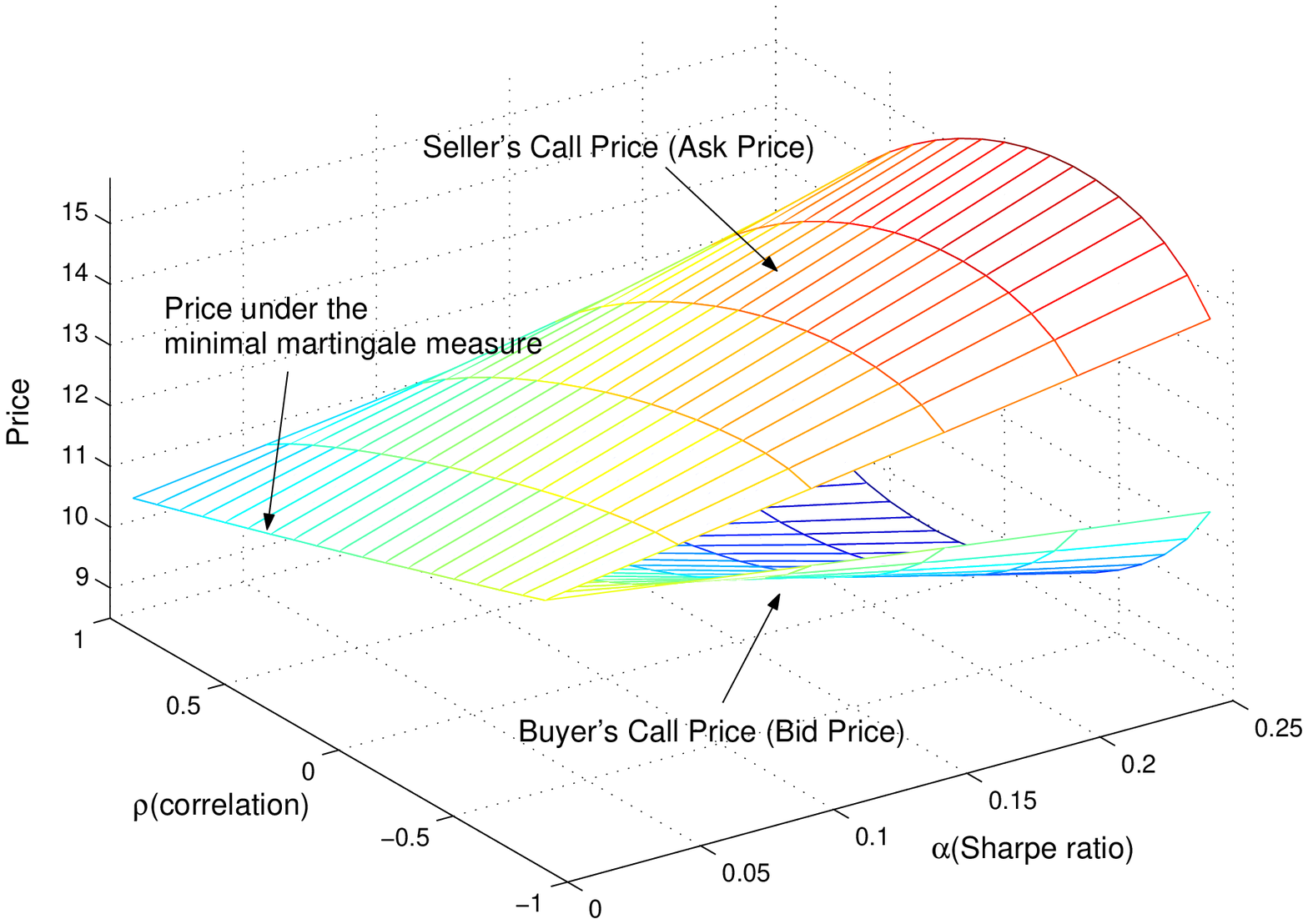,height=2.5in,width=3.5in} \psfig{figure=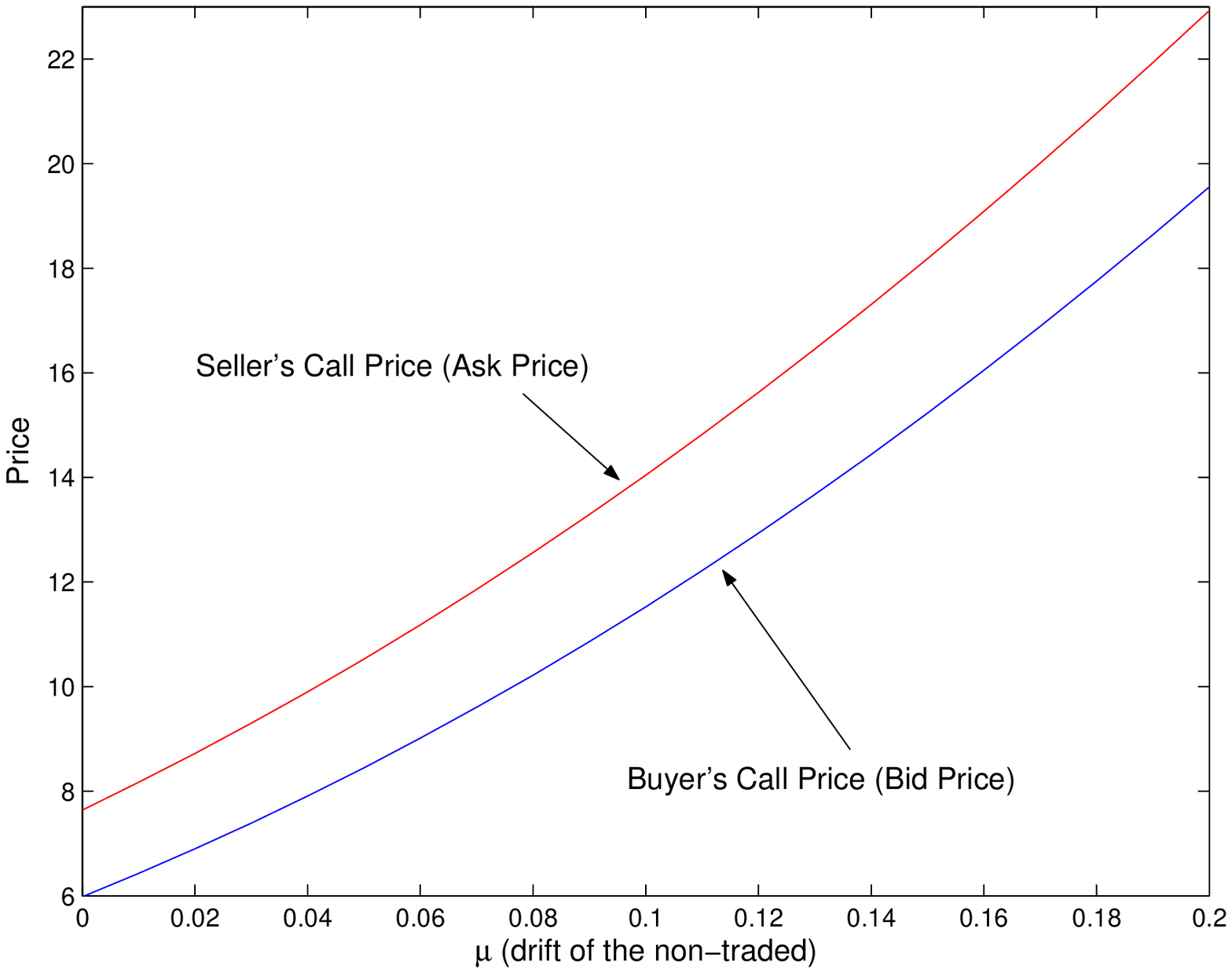,height=2.3in,width=2in}}
   \hfill \break
    Figure 1: Variation of the Bid and Ask Prices of a European Call Option with respect to $\rho$, $\alpha$ and $\mu$.}}

\noindent {\bf Remark:}  In this remark, we show how the bid/ask prices are related to the good deal bounds of Cochrane and Sa\'{a}-Requejo (2000).  To this end, note that (2.14) can be written as
$$
\left\{ \eqalign{&P_t + \left( \mu(S, t) - (a(H, t) - r) \rho {\sigma(S, t) \over b(H, t)} \right) S P_S + \max_{0 \leq |h|\leq \alpha} \left\{h \, \sqrt{1 - \rho^2} \, \sigma(S, t) \, S P_{S}\right\} + r H P_H \cr
& \quad + {1 \over 2} \sigma^2(S, t) S^2 P_{SS} + \rho \sigma(S, t) b(H, t) S H P_{SH} + {1 \over 2} b^2(H, t) H^2 P_{HH} - rP =0, \cr
& P(S, H, T) = g(S).} \right.$$

\noindent To obtain (2.30) for the buyer's price, replace $\max$ by $\min$ and replace $\a$ by $\beta$ in this Hamilton-Jacobi-Bellman equation. Therefore, $P$ and $P^b$ can be represented as
$$
P(S,H,t)=\max_{|h|\leq \alpha}{\bf E}^Q\left(e^{-r(T-t)}g(S_T)\bigg|{\cal F}_t\right),\quad P^b(S,H,t)=\min_{|h|\leq \beta}{\bf E}^Q\left(e^{-r(T-t)}g(S_T)\bigg|{\cal F}_t\right),
$$
in which $h = \{h_s\}_{0\le s \le T}$, and $|h|\leq \alpha$ means that we constrain $h$ so that $|h_s| \le \a$ for all $0 \le s \le T$. Here, the dynamics of the underlying processes under the equivalent martingale measure ${\bf Q}$ are
$$
\eqalign{dS_{t} &= \left(\mu(S_t,t)-(a(H_t,t)-r)\rho {\sigma(S,t) \over b(H_t,t)}  +h_t\sqrt{1-\rho^2} \, \sigma(S,t)\right)S_t \, dt \cr
& \quad + \sigma(S_t,t) \, S_t \left(\rho dB_t+ \sqrt{1-\rho^2}dB^{\perp}_t \right),}
$$
and
$$
dH_t = r \, H_t \, dt + b(H_t,t) \, dB_t,
$$
where $B$ and $B^{\perp}$ are independent Brownian motions under ${\bf Q}$, $B_t = W_t + \int_0^t {(a(H_s) - r) \over b(H_s, s)} \, ds$, and $B^{\perp}_t = W^{\perp}_t - \int_0^t h_s \, ds$.

The price $P$ is arbitrage free since it lies in the no-arbitrage price interval
$$
\left(\inf_{Q \in {\cal M}}{\bf E}^{Q}\left(e^{-r(T-t)}g(S_t)\bigg| {\cal F}_t\right), \sup_{Q \in {\cal M}}{\bf E}^{Q}\left(e^{-r(T-t)}g(S_t)\bigg| {\cal F}_t\right)\right),
$$
in which ${\cal M}$ is the set of equivalent martingale measures; see Schachermayer (2000). Since the no-arbitrage pricing interval is too wide and practically useless, it was Cochrane and Sa\'{a}-Requejo (2000)'s idea to find a subinterval that is tight enough.  They chose the subinterval to be
$$
\left(\min_{|h|\leq \beta}{\bf E}^Q\left(e^{-r(T-t)}g(S_T)\bigg|{\cal F}_t\right), \;  \max_{|h|\leq \alpha}{\bf E}^Q\left(e^{-r(T-t)}g(S_T)\bigg|{\cal F}_t\right)\right) \eqno(2.33)
$$
Although these ideas are not precisely stated in Cochrane and Sa\'{a}-Requejo (2000), this is indeed how they construct a no-arbitrage sub-interval; see Bj\"{o}rk and Slinko (2006). Cochrane and Sa\'{a}-Requejo (2000) initially aimed to find upper and lower bounds for prices by putting a constraint on the Sharpe ratio (by ruling out ridiculously good deals). However, this formulation of the problem is intractable, and they ended up constructing the sub-interval above, by imposing a constraint on $h$, the volatility/market price of risk with respect to $B^{\perp}$ of the Radon-Nikodym derivative of the measure ${\bf Q}$ with respect to ${\bf P}$.

We independently arrive at the same sub-interval in (2.33) in Sections 2.2 and 2.5. The lower end point and the upper end point of the above subinterval are the buyer's price and seller's price, respectively, when these investors' risk preferences are characterized by the instantaneous Sharpe ratio of a suitably defined portfolio. Our contribution over Cochrane and Sa\'{a}-Requejo (2000) can be summarized as follows: (1) By deriving the lower and upper boundaries independently and by a completely different mechanism, we provide a meaning to the these terms, which are otherwise technical, as buyer (bid) and seller (ask) prices, and interpret $\a$ as a risk aversion parameter.  (2) We give conditions under which the PDE that one expects these boundaries to satisfy have unique solutions.  (3) We analytically show how the price changes with respect to the the model parameters.  (4) We show that the price under the minimal martingale measure always lies between the lower and upper bounds of Cochrane and Sa\'{a}-Requejo (2000) for any $\a$.  (5) In the next section, in a stochastic volatility framework, we show that the instantaneous Sharpe ratio is equivalent to the market price of volatility risk when the payoff of an option on the underlying is convex. This implies that the upper and lower bounds in (2.33) are attained by a martingale measures ${\bf Q}_0$ and ${\bf Q}_1$ whose Radon-Nikodym derivatives with respect to ${\bf P}$ have market price of volatility risk equal to $\a$ and $-\beta$, respectively.

\sect{3. Pricing Derivative Securities in the Presence of Stochastic Volatility}

In this section, we model the price $S$ of the underlying asset as a diffusion process solving
$$
\eqalign{dS_t &= \mu \, S_t \, dt + \b(\sigma_t) \, S_t \, dB_t, \cr
d \sigma_t &= a(\sigma_t,t) \, dt+ b(\sigma_t,t) \, dW_t,}
\eqno(3.1)$$

\noindent in which $B$ and $W$ are standard Brownian motions on a probability space $(\Omega, {\cal F}, {\bf P})$, which are correlated with correlation coefficient $\rho$, and $\b(\cdot) > 0$ is non-decreasing.  We will see that the instantaneous Sharpe ratio is equivalent to the market price of volatility risk, which we establish by showing that the nonlinear PDE in (3.8) below can be reduced to the linear PDE in (3.9). We assume that the coefficients $\beta$, $a$, and $b$ satisfy growth conditions and are locally Lipschitz in $\sigma$.  These conditions ensure the existence and uniqueness of solutions of the two SDEs in (3.1).

As in the previous section, suppose that a portfolio $\Pi$ contains an option liability with payoff $g(S_T)$ and a self-financing portfolio of $\pi_t$ shares of the risky asset and shares of a money market account that earns at the rate $r \ge 0$. Let $P(S,\sigma,t)$ be the price of the option at time $t$ when the stock price is equal to $S$ and the volatility is $\b(\sigma)$. By It\^o's Lemma, the value of a portfolio at time $t+h$ in terms of its value at time $t$ equals
$$
\eqalign{\Pi_{t+h} &= \Pi_t - \int_t^{t+h} {\cal D}^\mu P (S_s, \sigma_s, s) ds + \int_t^{t+h} \b(\sigma_s)S_s (\pi_s- P_S(S_s,\sigma_s,s) dB_s \cr
& \quad - \int_{t}^{t+h} b (\sigma_s, s) P_{\sigma} (S_s,\sigma_s,s) dW_s + \int_{t}^{t+h} [(\mu-r) \pi_s S_s+r \Pi_s] ds,} \eqno(3.2)$$

\noindent where the differential operator ${\cal D}^m$, for some constant $m \in \bf R$, is given by
$$
{\cal D}^m  \nu(S, \sigma, t) = \nu_{t}+ m S \nu_s + a(\sigma,t) \nu_{\sigma} + {1 \over 2} \b(\sigma)^2 S^2 \nu_{SS} + \rho \beta(\sigma) b(\sigma, t) S \nu_{S \sigma} + {1 \over 2} b^2(\sigma,t) \nu_{\sigma \sigma}-r \nu, \eqno(3.3)$$

\noindent  The value of $\pi_t$ that minimizes the local variance is given by
$$
\pi^*_t =  P_S (S_t, \sigma_t, t) + \rho \, {b(\sigma_t,t) \over \beta(\sigma_t) S_t} \, P_{\sigma}(S_t,\sigma_t,t), \eqno(3.4)$$

\noindent which is delta hedging plus a term that takes into account the volatility risk. Under this best local hedge, the portfolio's drift and local variance become
$$
\lim_{h \rightarrow 0} {1 \over h} ({\bf E}(\Pi_{t+h}|{\cal F}_t) - \Pi) = - {\cal D}^r P(S, \sigma, t) + \rho \, (\mu - r) \, {b(\sigma, t) \over
\beta(\sigma)} \, P_{\sigma}(S, \sigma, t)+r \Pi, \eqno(3.5)
$$
and
$$
\lim_{h \rightarrow 0} {1 \over h} {\bf Var}(\Pi_{t+h}|{\cal F}_t)= (1-\rho^2) b^2(\sigma, t) P^2_{\sigma}(S, \sigma, t). \eqno(3.6)$$

Because the minimum variance is positive, the writer of the option is not able to hedge the risk completely. Therefore, the option should be priced in such a way that the writer of the option is rewarded for taking the extra risk. Here, the price will depend on the risk preference of the option writer via the instantaneous Sharpe ratio.  As in Section 2, the price is determined from
$$
-{\cal D}^r P(S, \sigma,t) + {(\mu-r) \rho b(\sigma,t) \over \beta(\sigma)} P_{\sigma}+r \Pi = r \Pi + \a \sqrt{1 - \rho^2} \, b(\sigma, t) \big| P_{\sigma}(S, \sigma, t) \big|,  \eqno(3.7)$$

\noindent where $\a$ is the instantaneous Sharpe ratio of the investor. It follows from this equation that the writer's price of the option solves the following non-linear PDE:
$$
\left\{\eqalign{ &P_t+ r S P_S + \left(a(\sigma,t)- {(\mu-r) \rho b(\sigma,t) \over \beta(\sigma)}\right)P_{\sigma} + \alpha
\sqrt{1-\rho^2} \, b(\sigma,t) |P_{\sigma}| \cr
& \quad + {1 \over 2} \beta(\sigma)^2 S^2 P_{SS} +\rho \beta(\sigma) b(\sigma,t) S P_{S \sigma} + {1 \over 2} b^{2}(\sigma,t) P_{\sigma \sigma}-rP =0 \cr
& P(S, \sigma, T) = g(S).} \right. \eqno(3.8)$$

\noindent  If the investor is allowed to trade in an option with larger maturity or a secondary asset that has both Brownian motions present in its dynamics (the option has this type of dynamics because of It\^{o}'s lemma) besides the underlying asset, then the hedging error (or local variance) will be equal to zero almost surely.  In this case, $\alpha$ does not matter and one reduces to the complete market case--in particular, the solution of (3.8) is the price determined by the complete market.  However, in this paper, we focus on the incomplete market case and do not assume that the option writer can trade in such an asset.

As in Section 2, we have an existence and uniqueness result, as well as a comparison theorem, that applies to the PDEs in (3.8) and in (3.9) below; see Appendix B for details.

\medskip

\noindent{\bf Assumption 3.1.}  In the rest of this section, unless otherwise stated, we assume that the payoff function $g$ is twice differentiable and convex.

\medskip

Consider the following linear PDE:
$$
\left\{\eqalign{ &G_t+ r S G_S + \left(a(\sigma,t)- {(\mu-r) \rho b(\sigma,t) \over \beta(\sigma)} + \alpha \sqrt{1-\rho^2} \, b(\sigma, t) \right)G_{\sigma} \cr
& \quad + {1 \over 2} \beta(\sigma)^2 S^2 G_{SS} +\rho \beta(\sigma) b(\sigma,t) S G_{S \sigma} + {1 \over 2} b^2(\sigma, t) G_{\sigma \sigma} - rG=0  \cr
&G(S, \sigma, T) = g(S).} \right. \eqno(3.9)$$

\noindent  In a series of results (Theorems 3.3 and 3.5 and Corollary 3.7), we show that the solution of (3.9) equals the solution of (3.8). For the remainder of this section, we also make the following assumptions about the coefficients of (3.9).

\medskip

\noindent{\bf Assumption 3.2.}
\item{(1)} The functions $a(\sigma,t)$, $b(\sigma,t)$, and $\beta(\sigma)$ are bounded in ${\bf R} \times [0,T]$ or $\bf R$, as appropriate.

\item{(2)} (Uniform ellipticity) There exists a constant $\delta > 0$ such that $\max(\beta (\sigma), b(\sigma,t)) \ge \delta$ for all $\sigma \in {\bf R}$ and $t \in [0,T]$.

\item{(3)} Let
$$
\gamma(\sigma,t) = a(\sigma,t)- {(\mu-r) \rho b(\sigma,t) \over \beta(\sigma)} + \alpha \sqrt{1-\rho^2} \, b(\sigma,t); \eqno(3.10)$$
then, $\beta$, $b$, and $\gamma$ are Lipschitz, uniformly in $t$.  Because we assume that $\beta$ is bounded from above and below and that $b$ is bounded from above, the Lipschitz assumption on $\gamma$ can be achieved by assuming that $a$ is Lipschitz, uniformly in $t$.

\item{(4)} The payoff $g$ satisfies the following growth condition: There exist constants $L > 0$ and $\la \ge 1$ such that $$g(S) \le L \left(1+ (\ln S)^{2 \lambda}\right) $$ for all $S > 0$.

\th{3.3}{ Under Assumption $3.2$, the PDE in $(3.9)$ has a unique $C^{2, 2,1}({\bf R^+} \times {\bf R} \times [0,T])$ solution in ${\cal C},$ which is defined in Theorem $B.1$. Moreover, its derivative satisfies the same growth condition as given in $(B.2)$.}

\pf The proof follows from Karatzas and Shreve (1991, pages 366-367) and Friedman (1975, Theorem 5.3). $\square$

\medskip

\noindent{\bf Assumption 3.4.} In the rest of the section, we assume the following statements hold:

\item{(1)} $\beta$, $b$ and $\gamma $ are differentiable with respect to $\sigma$.

\item{(2)} The derivatives $\beta_{\sigma}$, $b_{\sigma}$ and $\gamma_{\sigma}$ are bounded, and $\beta_\sigma$ and $b_{\sigma}$ are Lipschitz with respect to $\sigma$ uniformly in $t$.

\item{(3)} The terminal payoff function $g$ is $C^{2}$, and there exist constants $M > 0$ and $\mu \ge 1$ such that
$$
g''(S) \le M \left(1+ ( \ln S)^{2 \mu}\right) $$
for all $S > 0$.

\th{3.5}{ Under Assumptions $3.1,$ $3.2,$ and $3.4,$ the solution $G$ of $(3.9)$ satisfies $G_\sigma > 0$ on ${\bf R}^+ \times {\bf R} \times [0, T]$.}

\pf By differentiating the PDE that $G$ satisfies with respect to $\sigma$, we obtain that $f = G_\sigma$ solves
$$
\left\{\eqalign{ & f_t + (r + \rho (\beta_{\sigma}(\sigma) b(\sigma,t)+ \beta(\sigma) b_{\sigma}(\sigma,t))) S f_S + (\gamma(\sigma,t)+ b(\sigma,t) b_{\sigma} (\sigma,t))f_{\sigma} \cr
& \qquad + {1 \over 2} \beta(\sigma)^2 S^2 f_{SS} + \rho \sigma S b(\sigma,t) f_{S \sigma} + {1 \over 2} b^2(\sigma, t) f_{\sigma \sigma}+ (\gamma_{\sigma}(\sigma,t)-r) f \cr
& \quad = -\beta(\sigma) \beta_{\sigma}(\sigma) S^2 G_{SS} \cr
& f(S, \sigma,T) = 0.} \right. \eqno(3.11)$$

Under Assumption 3.4(1) and (2), we can apply the Feynman-Kac Theorem (Karatzas and Shreve, 1991, Theorem 5.7.6) to represent the solution of this linear PDE as an expectation as follows:
$$
f(S,\sigma,t)= {\bf \tilde E}^{S,\sigma,t}\left[ \int_t^T \b(\sigma_s) \b_\sigma(\sigma_s) S^2_s G_{SS}(S_s,\sigma_s,s)
\exp\left(\int_t^s (\gamma_{\sigma}(\sigma_u, u )-r)du\right)ds\right], \eqno(3.12)$$

\noindent where under the measure $\tilde{\bf P}$ the stock price process is given by
$$
\eqalign{ &dS_t= (r + \rho (\beta_{\sigma}(\sigma_t) b(\sigma_t,t)+ \beta(\sigma_t) b_{\sigma}(\sigma_t,t))) S_t dt +
\b(\sigma_t) S_t d\tilde{B}_t, \cr
& d\sigma_t = (\g(\sigma_t, t)+ b(\sigma_t, t) b_{\sigma} (\sigma_t, t)) dt + b(\sigma_t,t) d\tilde W_t,} \eqno(3.13)$$

\noindent where $\tilde{B}$ and $\tilde{W}$ are standard Brownian motions with respect to $\tilde{\bf P}$ with correlation $\rho$.
Observe from this equation that if we can show that $G_{SS} > 0$, then we can conclude that $G_{\sigma} > 0$.

Similarly, by differentiating the PDE that $G$ satisfies twice with respect to $S$ and by using the Feynman-Kac representation,
we obtain that under a suitable probability measure $\hat{\bf P}$, the function $k = G_{SS}$ is given by
$$k(S, \sigma, t) = \hat{\bf E}^{S,\sigma,t} \left[ g''(S_T) \exp\left(\int_{t}^{T}(r+\beta^2(\sigma_u))du\right) \right] \ge 0. \eqno(3.14)$$
Note that Assumption 3.4(3) is only used in this last step.  $\square$

\cor{3.6} {Under Assumptions $3.1,$ $3.2,$ and $3.4$ and under the condition in Theorem $B.1[v]$ on the payoff function $g,$ the
solution $P$ of $(3.8)$ equals the solution $G$ of $(3.9)$.}

\pf Under the assumptions of the corollary, the assumptions of Theorem B.1 are satisfied; therefore, there is a unique
viscosity solution of this PDE satisfying a polynomial growth condition.  By Theorem 3.5, we have $G_\sigma > 0$, which implies
that $G$ solves (3.8). Because $G$ also satisfies a polynomial growth condition, it is the unique solution of (3.8).  $\square$

\medskip

\noindent{\bf Remark:} If $g' \ge 0$ on ${\bf R}^+$, then $P_S \ge 0$ by an argument similar to the derivation of (3.14) and $P_\sigma > 0$ by Theorem 3.5 and Corollary 3.6.  Then, the optimal investment in the risky asset $\pi^* > 0$ from (3.4); that is, there is no short selling in this case.

\medskip

\noindent {\bf Remark:}
\item{(1)} By using the Feynman-Kac formula, we obtain that the price $P$ is given by the expectation
$$
P(t,S,\sigma)={\bf \bar{E}}^{S, \sigma, t} \left[ e^{-r(T-t)} g(S_T) \right],$$
where under the measure ${\bf \bar P}$ the stock price process is given by
$$
\eqalign{dS_t &= r S_t + \beta(\sigma_t) S_t d \bar B_t, \cr
d \sigma_t &= \left(a(\sigma_t,t)- {(\mu-r) \rho b(\sigma_t,t) \over \beta(\sigma_t)} + \alpha \sqrt{1-\rho^2} \, b(\sigma_t, t) \right) dt+ b(\sigma_t,t)d \bar W_t,} \eqno(3.15)$$
where $\bar{B}$ and $\bar{W}$ are standard Brownian motions with respect to $\bar{\bf P}$ with correlation $\rho$.  The Radon-Nikodym derivative of $\bar{\bf P}$ with respect to ${\bf P}$ is given  by
$$
{d {\bf \bar{P}} \over {d \bf P}} \bigg|_{{\cal F}_t} = \exp \left( - \int_{0}^{t} {\mu-r \over \beta(\sigma_s)} dB_s -{1 \over 2} \int_{0}^{t} \left({\mu-r \over \beta(\sigma_s)}\right)^2 ds \right) \exp\left( \alpha W_t + {1 \over 2} \alpha^2 t \right).
\eqno(3.16)$$

\item{} Here $-\alpha$ is the market price of volatility risk. (See Fouque, Papanicolaou, and Sircar (2000, page 47) for the definition of the market price of volatility risk.) The price we obtain for an investor with a given risk tolerance, measured by the instantaneous Sharpe ratio, $\alpha$, coincides with the price under the risk-neutral measure with a market price of volatility risk that is equal to $-\alpha$.

\item{(2)} The results of this section extend to the case for which $\alpha = \alpha (\sigma,t)$ such that the function $\alpha$ is bounded and Lipschitz uniformly in $t$ and its derivative $\alpha_{\sigma}$ exists, is bounded, and is Lipschitz uniformly in $t$.

\medskip

\noindent {\bf Remark:}  We used Assumption 3.2(4) to establish the existence and uniqueness of a solution to the PDE in (3.9), which is a sufficient condition arising from the theory of parabolic PDEs (see Friedman (1975)). On the other hand, Assumption 3.4(3) is used in the maximum principle for the PDE of $G_{SS}$ (the second derivative with respect to $S$ of the solution of (3.9)).  Again, this condition comes from the theory of parabolic PDEs. However, Assumptions 3.2(4) and 3.4(3) are sufficient but not necessary for Corollary 3.6 (and also Theorem 3.5) to hold.  Indeed,

(i) Assumption 3.4(3) is not satisfied by the put option, but below in Theorem 3.7, we show that the conclusions of Corollary 3.6 and Theorem 3.5 hold for a European put option since the payoff function can be uniformly approximated by a sequence of smooth convex functions.

(ii) Neither Assumption 3.2(4) nor Assumption 3.4(3) holds for the call option.  But, by using put-call parity in Corollary 3.8 below, we show that the result of Corollary 3.6 is still valid.  Theorems 3.3 and 3.5 also hold.

\medskip

\th{3.7}{Suppose the assumptions of Corollary $3.6$ hold, except that $g(S)=(K-S)_+$.  Then, the result of Corollary $3.6$ still holds.}

\pf Let $(g_n)_{n \geq 0}$ be a sequence of bounded, convex, $C^2$ functions satisfying the assumptions of Corollary 3.6 that uniformly converge to $g(S)=(K-S)_+$. Then, $P^{n}(S, \sigma, t)={\bf \bar{E}}^{S, \sigma, t}\left[ e^{-r(T-t)} g_n(S_T)\right]$ converges uniformly to $P(S, \sigma, t)={\bf \bar{E}}^{S, \sigma, t}\left[ e^{-r(T-t)} g(S_T)\right]$. Since $(P^{n})_{n \geq 0}$ is a convex sequence of functions, as a result of (3.14), it follows that $P$ is also convex. Since $P$ is the unique smooth solution of (3.9) with $g(S)=(K-S)_+$ by Theorem 3.3, it follows that $P_{SS} \geq 0$.  Then, one can repeat the proofs of Theorem 3.5 and Corollary 3.6 without the condition in Assumption 3.4(3). $\square$

\cor{3.8}{Suppose the assumptions of Corollary $3.6$ hold, except that $g(S)=(S-K)_+$.  Then, the result of Corollary $3.6$ still holds.}

\pf Let $C(S, \sigma, t)={\bf \bar{E}}^{S, \sigma, t}\left[e^{-r(T-t)}(S_T-K)_+\right]$. Then, $C$ and $P$ from the proof of Theorem 3.7 satisfy the following parity relation: $C(S, \sigma, t) + K e^{-r(T - t)} = P(S, \sigma, t) + S$. Therefore, $C$ is the unique solution of (3.9) with $g(S) = (S - K)_+$, although it does not satisfy Assumption 3.2(4).  Moreover, from the parity between $C$ and $P$, it follows that $C_{SS}(S, \sigma, t) = P_{SS}(S, \sigma, t)>0$.  Then, once more, one can repeat the proofs of Theorem 3.5 and Corollary 3.6 without the condition in Assumption 3.4(3).  $\square$

\sect{4. Summary and Conclusions}

We expanded on a method to value risk in an incomplete market first introduced by Windcliff et al.\ (2007) and further developed by Young (2007).  We assume that the risk is ``priced'' via the instantaneous Sharpe ratio. Because the markets in which we price are incomplete, there is no unique pricing mechanism and one must assume something about how risk is valued.  One could use the principle of equivalent utility (see Zariphopoulou (2001) for a review) or the Esscher transform (Gerber and Shiu, 1994) to price the risk.

In this paper, we first applied our method to price options on non-traded assets for which there is a traded asset that is correlated to the non-traded asset. Our main contribution to this particular problem was to show that our seller/buyer prices are the upper/lower good deal bounds of Cochrane and Sa\'{a}-Requejo (2000) and of Bj\"{o}rk and Slinko (2006) and to analyze these prices.  Second, we applied our method to price options in the presence of stochastic volatility. Our main contribution to this problem was to show that the instantaneous Sharpe ratio is the negative of the market price of volatility risk, as defined in Fouque, Papanicolaou, and Sircar (2000).  In general, our pricing technique yields the good deal bounds of Cochrane and Sa\'{a}-Requejo (2000); thereby, we provided a different motivation for good deal bounds.

\sect{Appendix A: Comparison Principle for the Results in Section 2}

In this appendix, we present a comparison principle from Barles et al.\ (2003) on which we rely extensively in Section 2.

\th{A.1} {Let $\cal G$ denote the set of functions defined on $G = {\bf R^+} \times {\bf R^+} \times [0, T]$ that are twice differentiable in their first and second variables and once in their third.  Define a differential operator $\cal L$ on $\cal G$ by
$$
{\cal L} v = v_t + {1 \over 2} \sigma^2(S, t) S^2 v_{SS} + \rho \sigma(S, t) b(H, t) S H v_{SH} + {1 \over 2} b^2(H, t) H^2 v_{HH} + h(S, H, t, v, v_S, v_H), \eqno({\rm A}.1)$$
for some function $h$. Assume that the operator $\cal L$ satisfies the following properties:}

\item{$[i]$} {\it There exists a constant $C_1>0$ such that $|\sigma(S_1,t)-\sigma(S_2,t)|$ $\le C_1 |\ln S_1 - \ln S_2|$ for all $S_1, S_2 > 0,$ and there exists a constant $C_2 > 0$ such that $\sigma(S, t) \le C_2 \sqrt{S (1 + \ln S)}$ for all $S > 0$ and $t \in [0, T]$.}

\item{$[ii]$} {\it There exists a constant $C_3 > 0$ such that $|b(H_1, t) - b(H_2, t)| \le C_3 | \ln H_1 - \ln H_2|$ for all $H_1, H_2 > 0,$ and there exists a constant $C_4 > 0$ such that $b(H, t) \le C_4 \sqrt{H (1+ \ln H)}$ for all $H > 0$ and $t \in [0, T]$.}

\item{$[iii]$} {\it There exist functions $d_1$ and $d_2$ satisfying
$$
\eqalign{ & 0 \le d_1(S, H, t) \le K S (1 + |\ln S| + |\ln H|), \cr
& 0 \le d_2 (S, H, t) \le K H (1 + |\ln S| + |\ln H|), } \eqno({\rm A}.2)$$
such that
$$
\left |h(S,H,t,v,p,z)-h(S,H,t,v,q,w)\right| \le d_1(S,H,t) |p-q|+ d_{2}(S,H,t)|z-w|. \eqno({\rm A}.3)$$
for all $S, H > 0$ and $t \in [0, T]$.}

\item{$[iv]$} {\it There exists a constant $m_1 > 0$ such that
$$
\left| h\left(S_1,H_1,t,v, {p \over S_1}, {q \over H_1}\right)- h\left(S_2,H_2,t,v, {p \over S_2}, {q \over H_2} \right)\right| \le m_1 k(p,q) \left[\left| \ln {S_1 \over S_2} \right|+ \left| \ln {H_1 \over H_2} \right|\right], \eqno({\rm A}.4)$$
for any $S_1, H_1, S_2, H_2 > 0$, $t \in [0,T]$, and $p, q \in \bf{R}$, in which  $k(p, q)=1+ \sqrt{p^2 + q^2}$.}

\item{$[v]$} {\it There exist constants $\gamma \in [0, (1 + \sqrt{5})/2)$ and $m_2 > 0$ such that
$$
|g(S_1)-g(S_2)|\le m_2 (1+ | \ln S_1 |+| \ln S_2 |)^{\gamma} | \ln S_1 - \ln S_2|,$$
for all $S_1, S_2 > 0$.}

\noindent {\it Denote by $\cal C$ the set of all locally bounded functions, $u$, that satisfy the following condition for some $k > 0$:
$$
{u(S, H, t)\over 1+ (|\ln S| + |\ln H|)^k} \rightarrow 0,$$

\noindent uniformly with respect to $t \in [0, T]$, as $|\ln S| + |\ln H| \rightarrow \infty$. Then, we can conclude the following two statements:

\noindent{$[$Existence and Uniqueness$]$} There exists a unique continuous viscosity solution in $\cal C$ of
 ${\cal L} v = 0$ with terminal condition $v(S,H,T) = g(S);$ see Crandall et al.\ $(1992)$ for the definition of a viscosity solution.

\noindent{$[$Comparison$]$} Let $u,v \in {\cal C}$ be continuous functions such that ${\cal L} u \ge 0 \ge {\cal L} v$ and $v(S, H, T) \le u(S, H, T)$ for all $S, H > 0$, then $v(S, H, t) \le u(S, H, t)$ for all $S, H > 0$ and $t \in [0,T]$.}

\medskip

\pf Transform the variables $S$, $H$, and $t$ in (A.1) to $x = \ln S$, $y = \ln H$ and $\tau = T - t$, and write $\tilde v(x, y, \tau) = v(S, H, t)$, etc.  Under this transformation, (A.1) becomes
$$
{\cal L} \tilde v = - \tilde v_\tau + {1 \over 2} \tilde \sigma^2(x, \tau) \tilde v_{xx} + \rho \tilde \sigma(x, \tau)
\tilde b(y, \tau) \tilde v_{xy} + {1 \over 2} \tilde b^2(y, \tau) \tilde v_{yy} + \tilde k(x, y, \tau, \tilde v, \tilde v_x, \tilde
v_y),  \eqno({\rm A}.5)$$

\noindent in which $\tilde k(x, y, \tau, \tilde v, \tilde p, \tilde z) = - {1 \over 2} \tilde \sigma^2(x, \tau) \tilde p - {1
\over 2} \tilde b^2(y, \tau) \tilde z + \tilde h(x, y, \tau, \tilde v, \tilde p, \tilde z)$, and $\tilde v$ is a differential
function on ${\bf R} \times {\bf R} \times [0, T]$.  Note that $P_S = e^{-x} \tilde P_x$ and $P_H = e^{-y} \tilde P_y$, so $p = e^{-x} \tilde p$ and $z = e^{-y} \tilde z$ in going from $h$ to $\tilde h$.  The differential operator in (A.5) is of the form considered by Barles et al.\ (2003); see that reference for the proof of our assertion.

The remaining item to consider is the form of the growth conditions in the original variables $S$, $H$, and $t$. Note that [i], [ii], and [v] are equivalent to assuming that $(x,t) \rightarrow \sigma(e^x,t)$ is Lipschitz in $x$ uniformly in $t$, $(x,t) \rightarrow b(e^x,t)$ is Lipschitz in $x$ uniformly in $t$, and the payoff function $x \rightarrow g(e^x)$ has polynomial growth.

From Barles et al.\ (2003), we know that analog of $(A.3)$ is
$$
\left |\tilde k(x, y, \tau, \tilde v, \tilde p, \tilde z) - \tilde k(x, y, \tau, \tilde v, \tilde q, \tilde w) \right| \le
\tilde d_1(x, y, \tau) \left|\tilde p -  \tilde q \right| + \tilde d_2(x, y, \tau) \left|\tilde z - \tilde w \right|, \eqno({\rm A}.6)$$

\noindent with
$$
\eqalign{& 0 \le \tilde d_1(x, y, \tau) \le K (1 + |x| + |y|), \cr
& 0 \le \tilde d_2(x, y, \tau) \le K (1 + |x| + |y|).} \eqno({\rm A}.7)$$

\noindent Under the original variables, the right-hand side of (A.6) becomes $d_1(S, H, t) |p - q| + d_2(S, H, t) \break |z - w|$, in which $d_1(S, H, t) =  \tilde d_1(x, y, \tau) e^x$ and $d_2(S, H, t) = \tilde d(x, y, \tau) e^y$ because $\tilde p = e^x p$ and $\tilde z = e^y z$.  Therefore, $\tilde d_1(x, y, \tau) \le K (1 + |x| + |y|)$ becomes $d_1(S, H, t) \le K e^x (1 + |x| + |y|) = K S (1 + | \ln S | + |\ln H|)$ and similarly for $\tilde d_2(x, y, \tau) \le K (1 + |x| + |y|)$. Note that the extra growth conditions in [i] and [ii] are to control the contribution from
the term $ - {1 \over 2} \tilde \sigma^2(x, \tau) \tilde p - {1 \over 2} \tilde b^2(y, \tau) \tilde z $ in the expression of
$\tilde k$.

Similarly, the condition in (A.4) can be drived from its analog in Barles et al.\ (2003), which is
$$
\left| \tilde k(x, y, \tau, \tilde v, \tilde p, \tilde z) - \tilde k(x', y', \tau, \tilde v, \tilde p, \tilde z) \right| \le
m_1 k(p,q) (|x - x'| + |y - y'|),$$

\noindent for some $m_1 > 0$. $\square$

\medskip

As a lemma for results in Section 2, we show that the differential operator associated with our problem satisfies the hypotheses of Theorem A.1.

\lem{A.2} {Define $h$ by
$$
h(S, H, t, v, p, z) = (\tilde{\mu}(S, H, t) + \a \sqrt{1 - \rho^2} \, \sigma(S, t) sgn(p)) S p + r H z - rv, \eqno({\rm A}.8)$$
\noindent in which $\tilde{\mu}(S, H, t) = \mu(S, t) - (a(H, t) - r) \rho {\sigma(S, t) \over b(H, t)}$.  Assume that
$\sigma$ satisfies Assumption $[i]$ in Theorem {\rm A.1}.  Then, $h$ satisfies Assumptions $[iii]$ and $[iv]$ of Theorem {\rm A.1} if we further assume that there exist constants $K, C_5 > 0$ such that
$$
|\tilde{\mu}(S, H, t) | \le K (1 + |\ln S| + |\ln H|),$$
and
$$
|\tilde{\mu}(S_1,H_1,t)-\tilde{\mu}(S_2,H_2,t)| \le C_5 | \ln S_1 - \ln S_2|, $$
for all $S_1, S_2, H_1, H_2, S, H > 0$ and $t \in [0, T]$.}

\pf
$$
\eqalign{& \big| h(S, H, t, v, p, z) - h(S, H, t, v, q, w) \big| = \big| \tilde{\mu}(S, H, t) S (p - q) + rH(z - w) \cr
& \qquad \qquad \qquad \qquad \qquad \qquad \qquad \qquad \qquad +  \a \sqrt{1 - \rho^2} \, \sigma(S, t) S (|p| - |q|) \big| \cr
& \quad \le \left( | \tilde{\mu}(S, H, t) | + \a \sqrt{1 - \rho^2} \, \sigma(S, t) \right) S |p-q| + rH |z - w|.} \eqno({\rm A}.9)$$

\noindent  Thus, (A.3) holds with $d_1(S, H, t) = ( | \tilde{\mu}(S, H, t) | + \a \sqrt{1 - \rho^2} \, \sigma(S, t)) S$ and $d_2(S, H, t) =  rH$.  Note that $d_2$ automatically satisfies (A.2), and $d_1$ satisfies (A.2) if $|\tilde{\mu}(S, H, t) | \le K (1 + |\ln S| + |\ln H|)$.

\noindent Next,
$$
\eqalign{  & \left| h\left(S_1, H_1, t, v, {p \over S_1}, {q \over H_1}\right)-  h\left(S_2,H_2, t, v, {p \over S_2}, {q \over
H_2}\right)\right| \cr
& \le \left[ \big| \tilde{\mu}(S_1, H_1, t) - \tilde{\mu}(S_2, H_2, t) \big|  + \alpha \sqrt{1-\rho^2} \, \big| \sigma(S_1, t) - \sigma(S_2, t) \big| \right] |p| .}$$

\noindent From this inequality, the assumption on $\sigma$ and the second assumption on $\tilde{\mu}$ it can be seen that (A.4) is satisfied with $m_1=C_5+C_{1} \alpha \sqrt{1-\rho^2}$.   $\square$

\sect{Appendix B: Comparison Principle for the Results in Section 3}

In this appendix, we present a comparison principle from Barles et al.\ (2003) on which we rely in Section 3.

\th{B.1} {Denote by $\cal L$ a differential operator on $\cal G$ whose action on a test function $v \in \cal G$ is given by
$$
{\cal L}v= v_{t}+ {1 \over 2} \beta(\sigma)^2 S^2 v_{SS} +\rho \beta(\sigma) S b  v_{S \sigma} + {1 \over 2} b^{2}(\sigma,t) v_{\sigma \sigma}+ h(S, \sigma, t , v, v_S, v_\sigma). \eqno({\rm B}.1)$$
\noindent We make the following assumptions about the operator $\cal{L}$.}

\item{$[i]$} {\it There exists a constant $c_1 > 0$, such that $|\beta(\sigma_1) - \beta(\sigma_2)| \le c_1 |\sigma_1 - \sigma_2|$ for any $\sigma_1, \sigma_2 \in {\bf R}$, and there exists a constant $c_2 > 0$ such that $|\beta(\sigma)| \le c_2 \sqrt{(1+|\sigma|)}$ for all $\sigma \in \bf R$.

\item{$[ii]$} There exists a constant $c_3 > 0$ such that $|b(\sigma_1,t) - b(\sigma_2,t)| \le c_3|\sigma_1 - \sigma_2|$ for any $\sigma_1, \sigma_2 \in {\bf R}$ and $t \in [0, T]$.

\item{$[iii]$} There exists a constant $c_4 > 0$ such that
$$|h(S, \sigma,t,v,p,z)-h(S,\sigma,t,v,q,w)| \le c_4 (1+ | \ln S|+ |\sigma|) (S |p-q| + |z-w|)$$
for any $S > 0,$ $\sigma, v, p, z, q, w \in {\bf R},$ and $t \in [0, T]$.

\item{$[iv]$} There exists a constant $c_5 > 0$ such that
$$\eqalign{&\left|h\left(S_1, \sigma_1, t, v , {p \over S_1}, q \right) - h\left(S_2, \sigma_2, t, v, {p \over S_2}, q \right)\right| \cr
& \qquad \le c_5 \left(1+ \sqrt{p^2+q^2} \right) \left(\left| \ln {S_1 \over S_2} \right| + |\sigma_1 - \sigma_2| \right)}$$
for any $S_1, S_2 > 0,$ $\sigma_1, \sigma_2, v, p, q \in \bf R,$ and $t \in [0, T]$.

\item{$[v]$} There exist constants $\gamma \in [0, (1 + \sqrt{5})/2)$ and $m_2 > 0$ such that
$$|g(S_1) - g(S_2)| \le m_2 (1+ | \ln S_1 |+| \ln S_2 |)^{\gamma} | \ln S_1 - \ln S_2|$$
for any $S_1, S_2 > 0$.

\noindent Denote by $\cal C$ the set of all locally bounded functions, $u$, that satisfy the following condition for some $k > 0$:
$${u(S, H, t)\over 1+ (|\ln S| + |\sigma|)^k} \rightarrow 0 \eqno({\rm B}.2)$$

\noindent uniformly with respect to $t \in [0, T]$, as $|\ln S| + |\sigma| \rightarrow \infty$.  Then, we can conclude the following two statements:

\noindent $[$Existence and Uniqueness$]$ There exists a unique continuous viscosity solution in $\cal C$ of ${\cal L}v = 0$ with terminal condition $v(S,\sigma,T)=g(S)$.

\noindent $[$Comparison$]$ Let $u,v \in {\cal C}$ be continuous functions such that ${\cal L} u \ge 0 \ge {\cal L} v$ and $v(S,
\sigma, T) \le u(S, \sigma, T)$ for all $S > 0, \sigma \in {\bf R}$, then $v(S,\sigma,t) \le u(S,\sigma,t)$ for all $S > 0, \sigma
\in \bf R,$ and $t \in [0,T]$.}

\medskip

\pf The proof follows from Barles et al.\ (2003, Theorem 2.1 and Corollary 2.1) after transforming the variables $S$ and $t$ in to
$x = \ln S$ and $\tau = T - t$ and is similar in nature to the proof of Theorem A.1; therefore, we omit the remainder of the details. $\square$

\medskip

We have the following lemma whose proof is similar to that of Lemma A.2, so we omit it.

\lem{B.2}{Define $h$ by
$$
h(S, \sigma, t, v, p, q)=  r S p + \sigma'(\sigma,t) q + \alpha \sqrt{1-\rho^2} \, b(\sigma,t) |q|-rv, \eqno({\rm B}.3)$$
in which $\sigma'(\sigma,t) = a(\sigma,t)- {(\mu-r) \rho b(\sigma,t) \over \beta(\sigma)}$. Then, $[iv]$ in Theorem $B.1$ is satisfied. Furthermore, if we assume that there exists a constant $K > 0$ such that
$$
|\sigma'(\sigma,t)| + \alpha \sqrt{1-\rho^2} b(\sigma,t)| \le K (1+|\sigma|), \eqno({\rm B}.4)$$
for any $\sigma \in {\bf R}$ and $t \in [0,T],$ then $[iii]$ in Theorem $B.1$ holds.}

\medskip

\centerline{\bf Acknowledgments} \medskip This research of the first author is supported in part by the National Science Foundation under grant DMS-0604491. We also thank the two anonymous referees for their very careful reading of our paper and their suggestions which improved our presentation.

\sect{References}

\smallskip \noindent \hangindent 15 pt Barles, G., S. Biton, M. Bourgoing, and O. Ley (2003), Uniqueness results for quasilinear parabolic equations through viscosity solutions' methods, {\it Calculus of Variations}, 18: 159-179.

\smallskip \noindent \hangindent 15 pt Bj\"{o}rk, T. and I. Slinko (2006), Towards a general theory of good deal bounds, {\it Review of Finance}, 10: 221-260.

\smallskip \noindent \hangindent 15 pt Cochrane, J. and J. Sa\'a-Requejo (2000), Beyond arbitrage: Good deal asset price bounds in incomplete markets, {\it Journal of Political Economy}, 108: 79-119.

\smallskip \noindent \hangindent 15 pt Cheridito, P. and M. Kupper (2006), Time-consistency of indifference prices and monetary utility functions, {\it preprint, Princeton University}, \hfill \break http://www.princeton.edu/$\sim$dito/papers/timf.pdf

\smallskip \noindent \hangindent 15 pt Crandall, M. G., H. Ishii, and P.-L. Lions (1992), User's guide to viscosity solutions of second order partial differential equations, {\it Bulletin of the American Mathematical Society}, 27 (1): 1-67.

\smallskip \noindent \hangindent 15 pt Davis, M. (2000), Optimal Hedging with Basis Risk, {\it working paper, Imperial College London}, http://www.ma.ic.ac.uk/~mdavis/docs/basisrisk.pdf

\smallskip \noindent \hangindent 15 pt F\"ollmer, H. and M. Schweizer (1991), Hedging of contingent claims under incomplete information, in: M. H. A. Davis and R. J. Elliott (eds.), {\it Applied Stochastic Analysis}, Stochastics Monographs, volume 5, Gordon and Breach, London/New York, 389-414.

\smallskip \noindent \hangindent 15 pt Forsyth, P. and G. Labahn (2006), Numerical methods for controlled Hamilton-Jacobi-Bellman PDEs in finance, {\it preprint, University of Waterloo}, \hfill \break
http://www.cs.uwaterloo.ca/~paforsyt/hjb.pdf.

\smallskip \noindent \hangindent 15 pt Fouque, J. P., G. Papanicolaou, and R. Sircar (2000), {\it Derivatives in Financial Markets with Stochastic Volatility}, Cambridge University Press, New York.

\smallskip \noindent \hangindent 15 pt Friedman, A. (1975), {\it Stochastic Differential Equations and Applications} 1, Academic Press, New York.

\smallskip \noindent \hangindent 15 pt Gerber, H. U. (1979), {\it Introduction to Mathematical Risk Theory}, Huebner Foundation Monograph 8, Wharton School of the University of Pennsylvania, Richard D. Irwin, Homewood, IL.

\smallskip \noindent \hangindent 15 pt Gerber, H. U. and E. S. W. Shiu (1994), Option pricing by Esscher transforms (with discussions), {\it Transactions of the Society of Actuaries}, 46: 99-191.

\smallskip \noindent \hangindent 15 pt Ilhan, A. and M. Jonsson, R. Sircar (2004), Portfolio optimization with derivatives and indifference pricing, to appear in {\it Volume on Indifference Pricing}, ed. R. Carmona, Princeton University Press; available at http://www.princeton.edu/$\sim$sircar/

\smallskip \noindent \hangindent 15 pt Karatzas, I. and S. E. Shreve (1991), {\it Brownian Motion and Stochastic Calculus}, second edition, Springer-Verlag, New York.

\smallskip \noindent \hangindent 15 pt Leung, T. and R. Sircar (2006), Accounting for risk aversion, vesting, job termination risk and multiple exercises in valuation of employee stock options, {\it preprint, Princeton University},
http://www.princeton.edu/$\sim$sircar/Public/ARTICLES/ESO\_utility.pdf.


\smallskip \noindent \hangindent 15 pt Musiela, M. and T. Zariphopoulou (2004), An example of indifference prices under exponential preferences, {\it Finance and Stochastics}, 8: 229-239.

\smallskip \noindent \hangindent 15 pt Royden, H. L. (1968), {\it Real Analysis}, second edition, Macmillan, New York.

\smallskip \noindent \hangindent 15 pt Schachermayer, W. (2000), Introduction to the mathematics of financial markets, {\it Lecture Notes in Mathematics}, no 1816, Springer Verlag, pp. 111-177.

\smallskip \noindent \hangindent 15 pt Schweizer, M. (2001), From actuarial to financial valuation principles, {\it Insurance: Mathematics and Economics}, 28: 31-47.

\smallskip \noindent \hangindent 15 pt Sircar R. and T. Zariphopoulou (2005), Bounds \& asymptotic approximations for utility prices when volatility is random, {\it SIAM Journal on Control and Optimization}, 43 (4): 1328-1353.

\smallskip \noindent \hangindent 15 pt Walter, W. (1970), {\it Differential and Integral Inequalities}, Springer-Verlag, New York.

\smallskip \noindent \hangindent 15 pt Windcliff, H., J. Wang, P. A. Forsyth, and K. R. Vetzal (2007), Hedging with a correlated asset: Solution of a nonlinear pricing PDE, {\it Journal of Computational and Applied Mathematics}, 200: 86-115.

\smallskip \noindent \hangindent 15 pt Young, V. R. (2007), Pricing life insurance under stochastic mortality via the instantaneous Sharpe ratio, working paper, Department of Mathematics, University of Michigan.

\smallskip \noindent \hangindent 15 pt Zariphopoulou, T. (2001), Stochastic control methods in asset pricing, in {\it Handbook of Stochastic Analysis and Applications}, D. Kannan and V. Lakshmikantham (editors), Marcel Dek-ker, New York.

 \bye